\documentclass{amsart}
\usepackage[foot]{amsaddr}
\usepackage{graphicx} % Required for inserting images
\usepackage{subcaption}
\usepackage{amsmath}
\usepackage{amsthm}
\usepackage{amsfonts}
\usepackage{xcolor}
\usepackage{comment}
\newcommand{\RR}{\mathcal{R}}

\DeclareMathOperator{\sign}{sgn}
\usepackage{epstopdf}
\AppendGraphicsExtensions{.eps}

\title[Self-similar imploding solutions]{Self-similar imploding solutions of the 1D compressible Euler equations with a far field cutoff}
\author{Jack Luong$^1$}
\thanks{Corresponding Author: Jack Luong}
\address{$^1$, Department of Mathematics, UCLA, Los Angeles, CA 90095, US}
\author{Scott Ramsey$^2$}
\address{$^2,^4$Los Alamos National Laboratory, Los Alamos, NM 87545, US}
\address{$^3$, Department of Mechanical Engineering, UCLA, Los Angeles, CA 90095, US}
\author{Andrea L. Bertozzi$^{1,3}$}
\author{Roy Baty$^2$}
\date{\today}
\email{jqluong@g.ucla.edu, ramsey@lanl.gov, bertozzi@g.ucla.edu, rbaty@lanl.gov}
\thanks{This work was supported by the US Department of Energy through the Los Alamos National Laboratory. Los Alamos National Laboratory is operated by Triad National Security, LLC, for the National Nuclear Security Administration of the US Department of Energy (Contract No. 89233218CN000001).}
\thanks{Jack Luong and Andrea L. Bertozzi are supported by Los Alamos National Lab award CW93449 C5677 and ONR grant N00014-23-1-2565.}
\thanks{A portion of this work was conducted while Andrea L. Bertozzi was an Ulam scholar at the Center for Nonlinear Studies at Los Alamos National Laboratory.}
 %Gas dynamics - general theory

\begin{document}
\begin{abstract}
Imploding solutions to the radially symmetric, isentropic, compressible Euler equations have been well-studied, inspired by the work of Guderley. However, these smooth imploding solutions are shown to be numerically unstable and difficult to compute in practice. On the other hand, the imploding solution of Kidder has a closed form solution and is numerically computable. But, it is unbounded in the far field. We consider Kidder's formulation in one dimension in which the unbounded far field condition is replaced with a constant density cutoff of the initial data. Strikingly, a non-centered rarefaction emerges from the cutoff and suppresses the implosion. We present an exact analytic solution to the problem with the cutoff and support our theoretical predictions with numerical simulations.
\end{abstract}
\maketitle
%Old abstract: Much recent investigation and progress has been made into smooth imploding solutions to the symmetric, isentropic, compressible Euler equations inspired by the work of Guderley. However, these smooth imploding solutions are shown to be numerically unstable and difficult to compute in practice. On the other hand, the imploding solution of Kidder has a closed form solution and is numerically computable, but is a nonphysical solution due to being unbounded radially. We propose modifying Kidder's solution in one dimension by continuously cutting off the density initial data by a constant after some point in space. From this initial condition, we construct the solution with rarefaction theory for systems of hyperbolic conservation laws. Strikingly, a rarefaction with nonstandard time dependence emerges and said rarefaction suppresses the implosion. We support our theoretical findings with numerical simulations.    
\section{Introduction}
Hyperbolic conservation laws are well known for \textit{finite time blowup} - the phenomenon where the solutions of hyperbolic conservation laws lose regularity in finite time. This normally takes the form of a \textit{shock}, as the gradient of solution variables becomes discontinuous. However, solutions to hyperbolic conservation laws may also blowup in finite time in the solution variables themselves. One such example are \textit{implosions} in fluid dynamics, where a flow directed towards the origin collapses into itself. Implosions appear in a wide range of contexts such as cavitation \cite{Hunter_cavity}, sonoluminescence \cite{Brenner_sonoluminescence}, gas cloud accretion of galaxies \cite{Woodward_implosion}, and star formation \cite{Hunter_supernova} along many others.

\par We are interested in implosions for the isentropic compressible Euler equations. The radially symmetric, compressible Euler equations detailing the conservation of mass, momentum, and entropy are given by
\begin{equation}\label{eq:con_mass}
    \frac{\partial \ln{\rho}}{\partial t} + v \frac{\partial \ln{\rho}}{\partial r} + \frac{\partial v}{\partial r} + \frac{nv}{r} = 0,
\end{equation}
\begin{equation}\label{eq:con_mom}
    \frac{\partial v}{\partial t} + v \frac{\partial v}{\partial r} + \frac{1}{\rho} \frac{\partial P}{\partial r} = 0,
\end{equation}    
\begin{equation}\label{eq:con_entropy}
    \frac{\partial}{\partial t} \ln{P\rho^{-\gamma}} + v\frac{\partial}{\partial r}\ln{P \rho^{-\gamma}} = 0,
\end{equation}
where $\gamma$ is the heat capacity ratio, $P$ is the pressure, $\rho$ is the density, and $v$ is the velocity. The spatial index $n$ is $0$ for planar geometry, $1$ for cylindrical geometry, and $2$ for spherical geometry. 

\par We restrict our attention to the case where entropy is initially constant in space. In doing so, we may eliminate the conservation of entropy equation.
Under these assumptions, the pressure is related to the density via $P(\rho) = C\rho^{\gamma}$ for some positive constant $C$. The compressible Euler equations under this assumption are the \textit{isentropic} compressible Euler equations. If the velocity is further assumed to be irrotational, the isentropic compressible Euler equations are also referred to as the \textit{potential flow} equations. In this work, since we study radially symmetric solutions to the compressible Euler equations, the solutions are automatically irrotational. The radially symmetric, compressible, isentropic Euler equations are
\begin{equation}\label{eq:isen_mass}
    \frac{\partial \ln{\rho}}{\partial t} + v \frac{\partial \ln{\rho}}{\partial r} + \frac{\partial v}{\partial r} + \frac{nv}{r} = 0,
\end{equation}
\begin{equation}\label{eq:isen_mom}
    \frac{\partial v}{\partial t} + v \frac{\partial v}{\partial r} + C\gamma\rho^{\gamma -2}\frac{\partial \rho}{\partial r} = 0.
\end{equation}
\par Since the flow is irrotational, we may define a scalar potential function $\phi$ such that $\nabla \phi = v$ and rewrite equations \eqref{eq:isen_mass} and \eqref{eq:isen_mom} with respect to the potential. We obtain the potential flow equations
\begin{equation}\label{eq:poten_mass}
    \rho_t + \nabla \cdot (\rho \nabla \phi) = 0
\end{equation}
and
\begin{equation}\label{eq:poten_mom}
    \phi_t + \frac{1}{2} |\nabla \phi|^2 = -C\frac{\gamma}{\gamma - 1} \rho^{\gamma - 1}. 
\end{equation}
Writing the isentropic, irrotational Euler equations in potential form gives us rich insight into the flow. We may interpret equation \eqref{eq:poten_mom} as a level set equation where level sets of the potential advect in their normal direction. However, the value of the level sets themselves change depending on the density. Equation \eqref{eq:poten_mass} dictates how density advects towards local maxima of the potential.   
\par Another physically relevant quantity is the local speed of sound, defined by 
\begin{equation}\label{eq:sos}
    c = \sqrt{\frac{\partial P}{\partial \rho}}.
\end{equation}
For potential flow, the speed of sound reduces to 
\begin{equation}\label{eq:sos_potential}
    c = \sqrt{\gamma}\rho^{\frac{\gamma - 1}{2}}.
\end{equation}

\subsection{Background}
\par A vast amount of literature on finite time blow ups exists for shocks for the Euler equations, and we discuss a few results though this discussion is by no means exhaustive. The first rigorous study of shocks for the Euler equations is attributed to Lax \cite{Lax1964} who used Riemann invariants to show the formation of shocks for small initial data of systems of two hyperbolic PDEs. Further developments in shock formation for the 1D isentropic Euler equations were advanced in \cite{Lebaud_thesis}. The authors of \cite{Chen_shock} extend the technique of Riemann invariants to show the formation of shocks in the isentropic Euler equations for large data.

\par Sideris \cite{Sideris_1985} was the first to show the formation of finite time singularities for the Euler equations in higher dimensions with constant far field initial conditions. He accomplished this with an energy argument. However, this argument cannot reveal the exact nature of the finite time singularity. In \cite{Yin_3DShock}, the author showed such a finite time singularity can be a shock for the 3D Euler equations with spherical symmetry. With the added assumption of irrotational solutions without symmetry neccessarily, the authors of \cite{Christodolou_irrotational} proved initial compression forms shocks for the Euler equations in higher dimensions with constant far field initial conditions.

\par On the other hand, less is known about implosions. A general theory for implosions for general hyperbolic equations has not been developed, and we instead outline the study of implosions for the compressible Euler equations. The modern study of implosions to the compressible Euler equations began with the seminal work of Guderley \cite{Guderley1942}. Guderley constructed preshocked initial data to the compressible Euler equations that forms an implosion at a later time. Moreover, the solution is also self-similar. By leveraging the self-similar variable, Guderley casts the problem into a system of ODEs. The corresponding phase portrait is used to construct the imploding solution. However, the initial data begins discontinuous, and it is unclear whether an implosion can arise from continuous and $C^1$ initial data - much like how a shock can develop from smooth initial data. 

\par The previous question has garnered recent interest building upon the self-similar reduction to a system of ODEs and phase portrait analysis framework of Guderley, particularly for the isentropic Euler equations. Jennsen and Tsikkou \cite{Jenssen2020} improve on the work of Guderley by constructing self-similar, radially symmetric, and continuous imploding solutions that have positive pressure everywhere, in contrast to the vanishing pressure fields constructed in \cite{Guderley1942}. However, the constructed solution has unbounded energy, although they state the energy lives in $L^1_{loc}$. The authors argue they can upgrade their solutions to finite mass and energy via a finite speed of propagation argument. They claim by choosing a proper modification of the initial data in the far field, the constructed solution near the origin will still persist and implode as desired at implosion time.

\par In the celebrated work of \cite{implosion_1} and \cite{implosion_2}, Merle et. al improve upon the imploding solutions of \cite{Guderley1942} and \cite{Jenssen2020} as these self-similar, radially symmetric, imploding solutions are both smooth and decay at infinity for the isentropic Euler and Navier-Stokes equations. The self-similar solutions constructed in \cite{implosion_1} are not of finite energy, but the authors are able to still construct finite energy, smooth, imploding solutions to the compressible Euler and Navier-Stokes Equations with the work done in \cite{implosion_2}. Merle et. al show this for almost every choice of heat capacity ratio $\gamma > 1$. However, one of the $\gamma$ omitted is the physically significant $\gamma = \frac{5}{3}$ - the heat capacity ratio of a monatomic gas. Buckmaster et. al \cite{buckmaster_implosion} improve on this work by constructing smooth imploding solutions for all $\gamma > 1$. They employ a computer-assisted proof to construct the imploding self-similar solutions and rigorously treat error bounds with interval arithmetic. This marks the first use of computer-assisted proof techniques to treat these types of problems. Building upon this work, the authors of \cite{caolabora2025nonradialimplosioncompressibleeuler} extend this methodology to construct self-similar solutions for the isentropic Euler and Navier-Stokes equations that are not radially symmetric on the torus and $\mathbb{R}^d$.

\par Although there have been many recent developments on implosions for the isentropic Euler equations, there still remain some practical concerns. For one, the solutions presented in the above works are not explicitly constructed. This leads to the second issue of the smooth, self-similar imploding solutions being computationally unstable as explored in Biasi's work \cite{biasi_compute}. Although Biasi details a numerical procedure on finding these smooth, self-similar imploding solutions, they tend to shock first when computed. These obfuscating details make it difficult to investigate these smooth imploding solutions physically or computationally. 

\par This motivates us to instead turn to the self-similar, radially symmetric, smooth, imploding solutions presented in the seminal work of Kidder in \cite{Kidder1974} and expanded upon in \cite{Kidder1976}. Motivated by the application of inertial containment fusion (ICF), Kidder in \cite{Kidder1974} investigates how isentropic compression can be used to severely heat and compress a pellet of thermonuclear fuel. Isentropic compression is desirable as it represents the maximum compression achievable for a given pressure. The gas is assumed to be inviscid and isentropic and thus modeled by the compressible, isentropic Euler equations. By considering the case of homogeneous isentropic compression - compression where each volume element is compressed identically - Kidder constructs an explicit, self-similar, radially symmetric, smooth, imploding solution to the compressible, isentropic Euler equations for $\gamma = \frac{5}{3}$. The concept of homogeneous compression originates with the work done in self-similar solutions in fluid mechanics as recorded by Sedov \cite{Sedov_book}. In Kidder's solution, the density and pressure approach infinity as time approaches the implosion time everywhere in the spatial domain except the origin, a distinct difference from the implosions occurring at a single point exhibited by the solutions presented in \cite{Guderley1942},\cite{Jenssen2020},\cite{implosion_1},\cite{buckmaster_implosion}, and \cite{caolabora2025nonradialimplosioncompressibleeuler}. Kidder expands his work to the case of a hollow shell pellet in \cite{Kidder1976}. These solutions have since enjoyed great use in ICF modeling, as cataloged in \cite{Meyer_book}. Both the Guderley and Kidder solutions have been viewed under the same framework by studying the dynamical system obtained via self-similar reduction as done in \cite{Meyer_1982}. 

\subsection{A potential flow perspective on implosions}
\par We note in general implosions are difficult to come by in the isentropic Euler equations. The prototypical example of an implosion is the unbounded accumulation of mass at a point as exhibited via the imploding Gaussian solution for the full Euler equations demonstrated in \cite{Ramsey_linear_velocity}. However, for the isentropic Euler equations, local maxima in density are discouraged from imploding. We can understand this by studying the potential flow equations \eqref{eq:poten_mass} and \eqref{eq:poten_mom} in one dimension. 

\par Suppose we have a stationary implosion at a local maximum for density located at $x^*$ in finite time. So, there exists a time $t^*$ and implosion time $T$ such that $\lim_{t \to T} \rho(x^*,t) = \infty$ with $u(x^*,t) = 0$ and $\rho(x^*,t)$ is a local maximum for $\rho$ for all times $t^* \leq t \leq T$. Studying equation \eqref{eq:poten_mass} at $x^*$ for times $t \geq t^*$ yields
\begin{equation}
    \rho_t(x^*,t) = -\rho(x^*,t) \phi_{xx}(x^*,t). 
\end{equation}
with solution
\begin{equation}
    \rho(x^*,t) = \rho(x^*,t^*) \exp{\int_{t^*}^t -\phi_{xx}(x^*,s) ds}.
\end{equation}

\par Now, in potential flow, gas flows in the direction of maxima of the potential. So, $\phi$ also has a local maximum at $x^*$ and the Laplacian of $\phi$ at $x^*$ is negative. Thus, the density can only continue to grow as long as $\phi_{xx}$ remains negative. For the density to implode in finite time, $\phi_{xx}$ itself must blow up in finite time.

\par The dynamics of $\phi_{xx}$ are obtained by taking the Laplacian of \eqref{eq:poten_mom}. We have
\[
    \partial_t\phi_{xx} + \phi_{xx}^2 + \phi_x \phi_{xxx} = -C\frac{\gamma}{\gamma - 1} \partial_{xx} \rho^{\gamma-1}
\]
which evaluated at $x^*$ is
\begin{equation}\label{eq:poten_mom_laplace}
    \partial_t\phi_{xx}(x^*,t) = -\phi_{xx}^2(x^*,t) -C\frac{\gamma}{\gamma - 1} \partial_{xx} (\rho(x^*,t)^{\gamma-1}).
\end{equation}
Examining the terms of the right hand side of equation \eqref{eq:poten_mom_laplace}, the Laplacian of $(\rho(x^*,t)^{\gamma-1})$ remains negative since $x^*$ is a local maximum of $\rho$ and $\gamma > 1$. In contrast, we know the term $-\phi_{xx}^2(x^*,t)$ causes finite time blowup left on its own. The growth of $\phi_{xx}(x^*,t)$ is balanced by the competing factors of $-\phi_{xx}^2(x^*,t)$ and $-\partial_{xx}(\rho(x^*,t)^{\gamma-1})$ with the pressure term counterbalancing the imploding term.   

\par
\begin{figure}[h]
     \centering
     \begin{subfigure}[h]{0.49\textwidth}
         \centering
         \includegraphics[width=\textwidth]{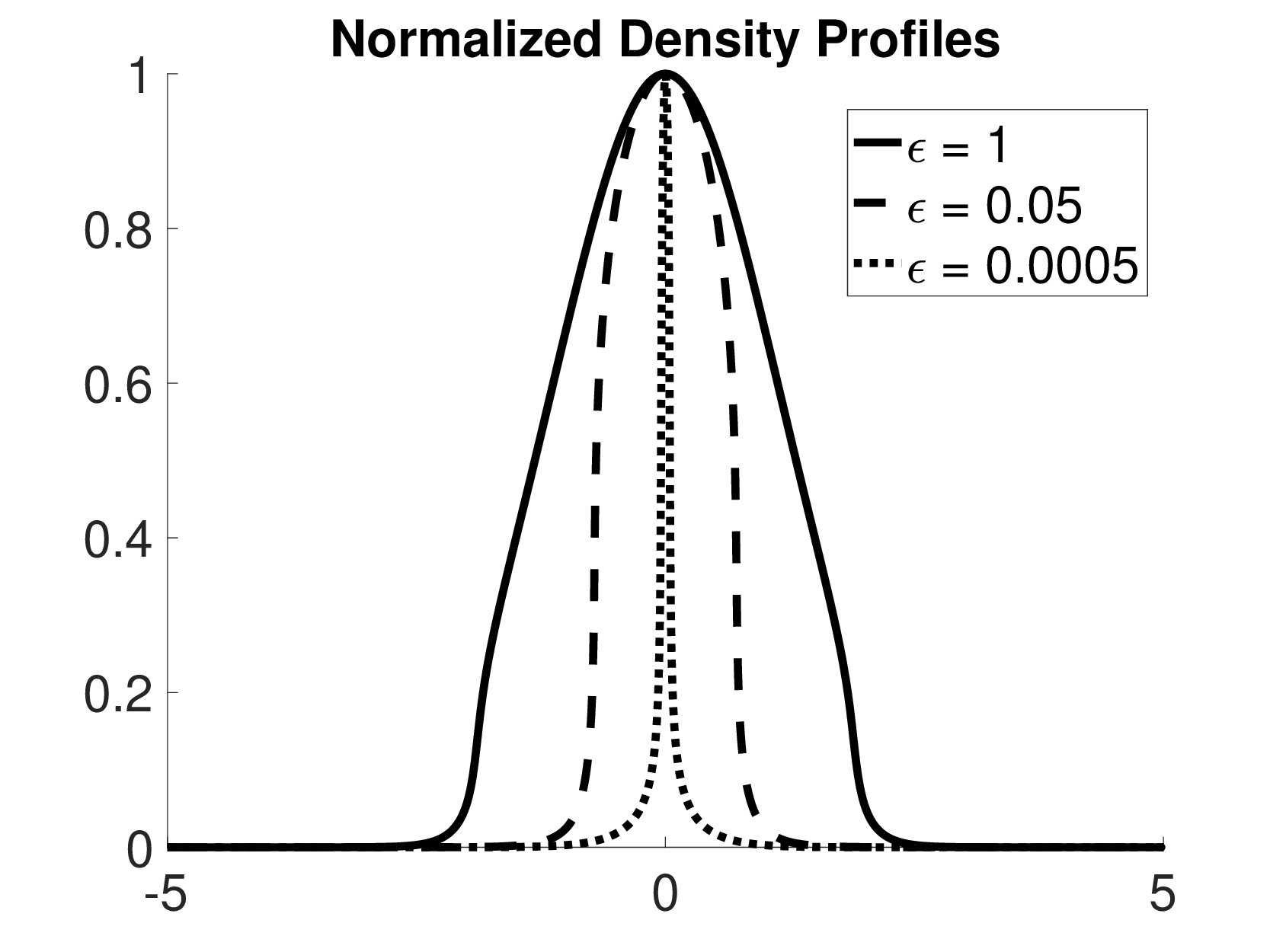}
         \caption{Density}
         \label{subfig:density_farfield}
     \end{subfigure}
     \hfill
     \begin{subfigure}[h]{0.49\textwidth}
         \centering
         \includegraphics[width=\textwidth]{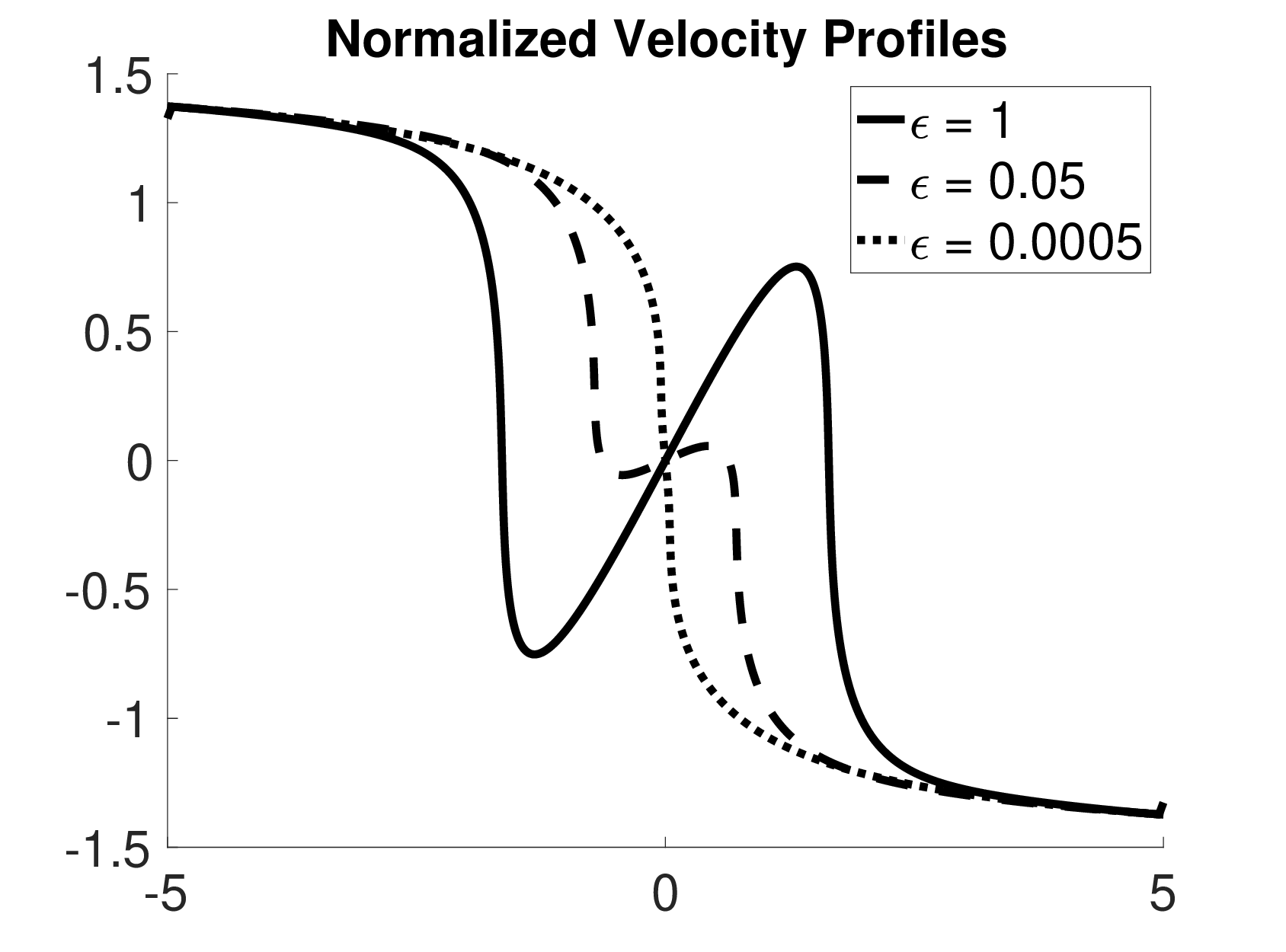}
         \caption{Velocity}
         \label{subfig:velocity_farfield}
    \end{subfigure}
    \caption{Normalized density and velocity for various choices of $\epsilon$ in initial data \eqref{eq:gaussian} at their respective shock times.}
    \label{fig:farfield}
\end{figure}

 For this reason, in practice it is difficult to observe these types of implosions. Consider the following initial data in 1D for equations \eqref{eq:isen_mass}  and \eqref{eq:isen_mom}:
\begin{equation}\label{eq:gaussian}
    \rho(x,0) = \epsilon \exp{(-x^2)}, \quad v(x,0) = \arctan(-x)
\end{equation}
where $\epsilon$ is a positive scalar. As $\epsilon \to 0$, equations \eqref{eq:isen_mass} and \eqref{eq:isen_mom} approach Burgers' equation. From numerical simulations, these solutions do not implode. For a short period of time, density starts to accumulate at the origin as desired, but after some critical time, the density starts to expand outwards from the origin and eventually collides with the density being directed towards the origin - forming a shock. Figure \ref{fig:farfield} shows the normalized density and velocity profiles for various choices of $\epsilon$ right before shock time. For any $\epsilon > 0$, two shocks form away at the origin representing the contact between gas expanding away from the origin and gas being compressed towards the origin. The shock position decreases as $\epsilon$ decreases. In general, as $\epsilon$ approaches $0$, the solutions approach the solutions of Burgers' equation. Additionally, the velocity profile begins to resemble the solution for Burgers' equation while the density approaches a delta function. 

\par
\begin{figure}[h]
     \centering
     \begin{subfigure}[h]{0.49\textwidth}
         \centering
         \includegraphics[width=\textwidth]{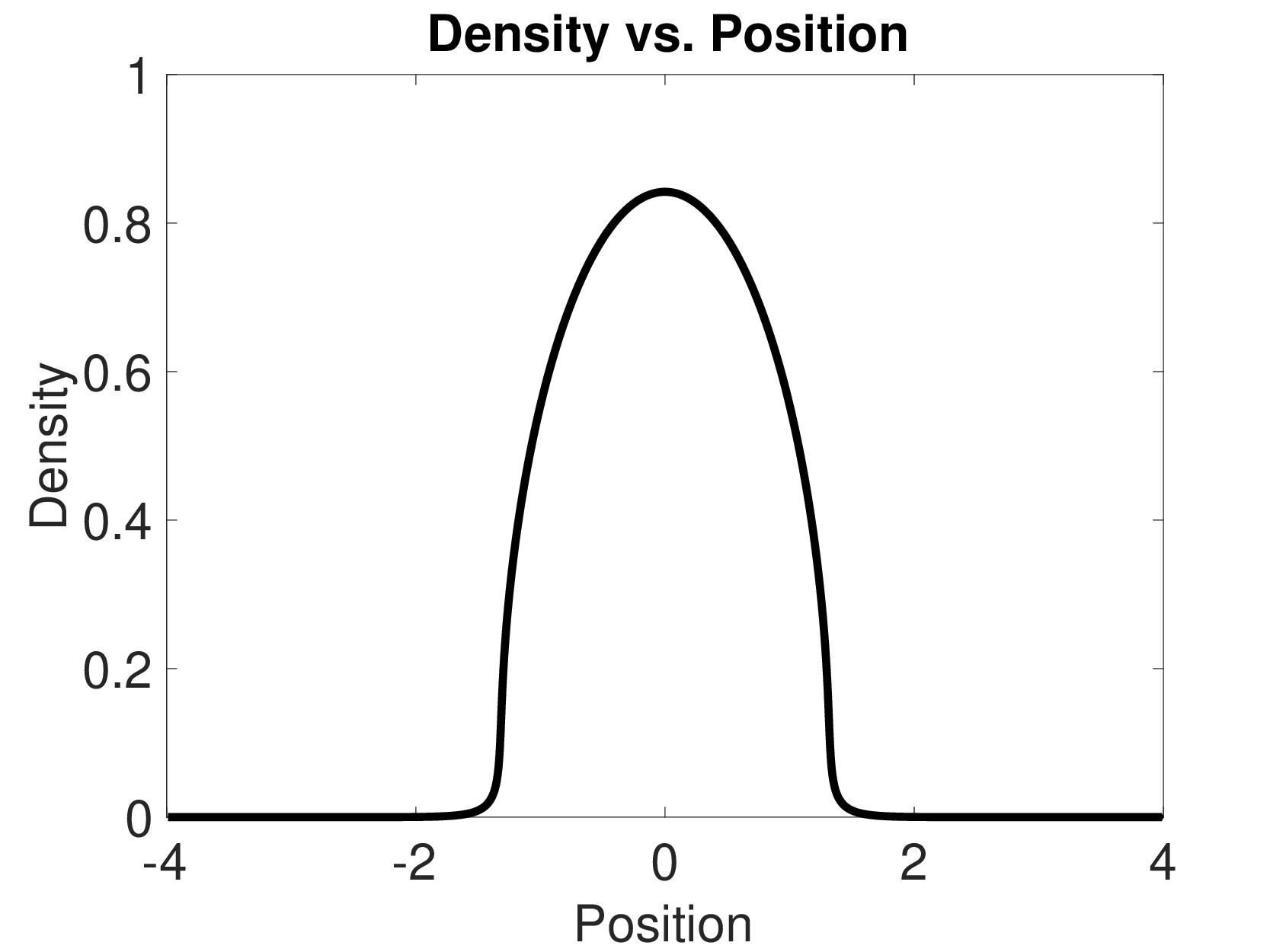}
         \caption{Density}
         \label{subfig:density_ps}
     \end{subfigure}
     \hfill
     \begin{subfigure}[h]{0.49\textwidth}
         \centering
         \includegraphics[width=\textwidth]{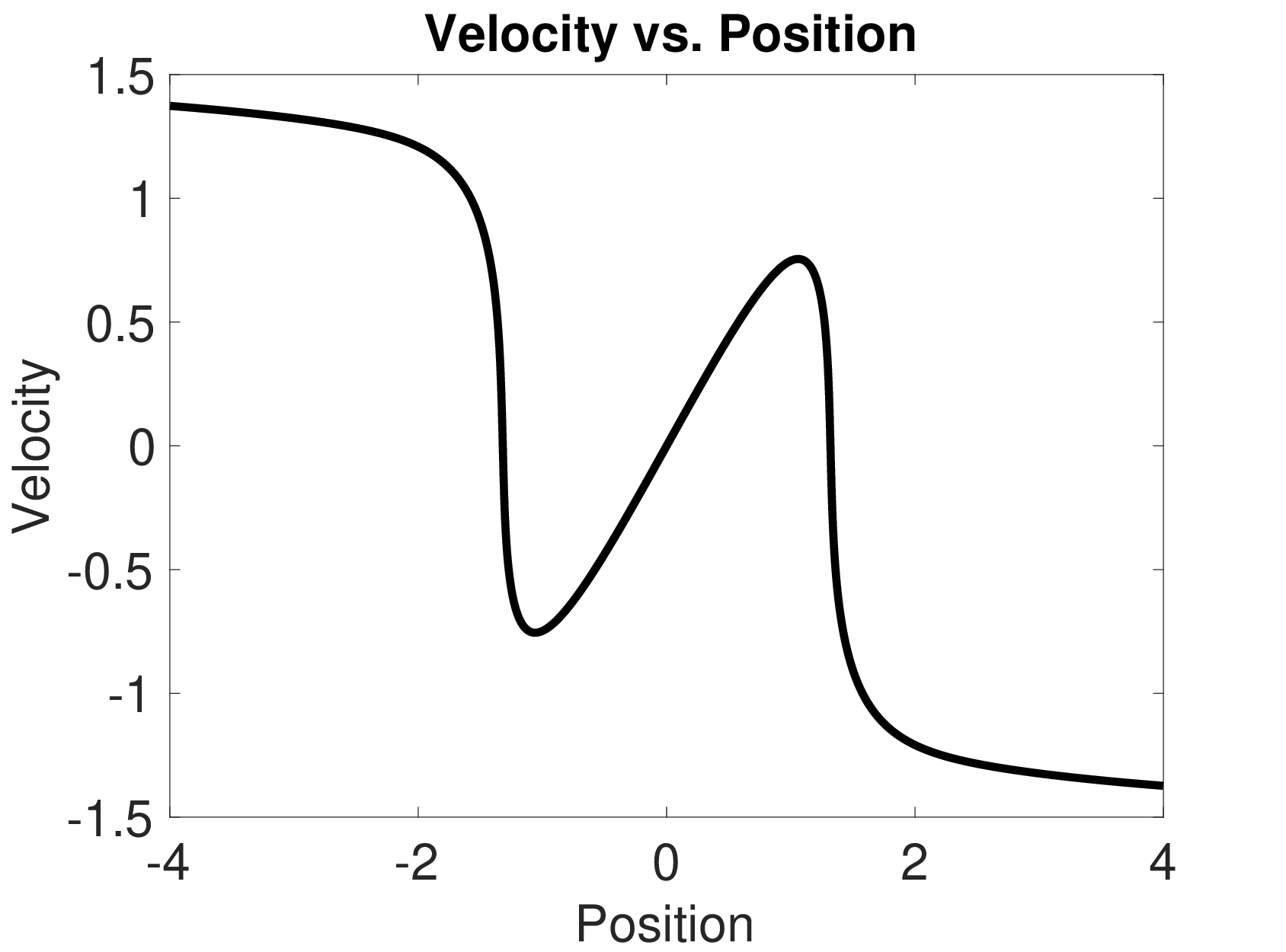}
         \caption{Velocity}
         \label{subfig:velocity_ps}
    \end{subfigure}
    \caption{Density and Velocity past the shock time for $\epsilon = 1$. After the shock occurs, the gas expands outwards into the compressing gas in the far field.}
    \label{fig:Postshock}
\end{figure}

 With the Gaussian density and arctangent velocity profile, we tried to induce an implosion by accumulating gas at the origin. The density has a local maximum at this accumulation point. However, the gas fails to implode, and a shock forms, indicating the gas expands away from the origin. Physically, the pressure of the gas at the origin grows until it grows too large and expands outwards. An example of this for $\epsilon = 1$ can be seen in Figure \ref{fig:Postshock}. The non-implosion occurs regardless of how much mass we start with. This numerical experiment, as well as the preceding discussion, showcases how constructing a simple imploding solution is difficult and naive approaches will not work. As far as the authors are aware, the literature studying this effect is lacking although numerical simulations are straightforward to perform. The challenges of forming an implosion via accumulation of mass at a stationary point along with the reported numerical instabilities by Biasi for the smooth imploding solutions of Merle et. al and Buckmaster et. al encourage further investigation of imploding solutions, such as the solution of Kidder. Kidder's solution avoids these pitfalls as the solution features an attractive global minimum in density rather than an attractive maximum.

\par Unfortunately, a major drawback to Kidder's solution is the solution is radially unbounded. Although it has an explicit form, the radial unboundedness makes it nonphysical. Since Kidder's imploding solution still has many desirable properties such as its explicit form and apparent numerical stability, we seek to modify the solution to make it more physically feasible while also preserving its desirable features. We propose continuously cutting off the initial data by a constant past some point in the domain. By doing so, we replace the nonphysical property of radial unboundedness with the milder physical requirement of enforcing constant density past a certain spatial point. Note this is the same setting studied in \cite{Sideris_1985} and \cite{Christodolou_irrotational}. The effect the constant cutoff has on the solution is governed by a finite speed of propagation. It is unclear how the cutoff will affect the imploding solution as it approaches implosion time and whether the implosion will still occur. This work seeks to answer such questions.

\par This paper is organized as follows. We begin by showing how the isentropic, homogeneous compression solution of \cite{Kidder1974} can be constructed via the principle of the linear velocity ansatz - a mathematical formalization of homogeneous compression - in Section \ref{sec:linear_vel_sol}. We state the general form of homogeneous isentropic compression solutions for arbitrary dimension and heat capacity ratio. Once the general form is obtained, we utilize the one dimensional solution as the basis for our subsequent analysis in Section \ref{sec:1D_cutoff}. We continuously cut off the initial data with a constant and derive the corresponding solution. To support our analysis, we present numerical solutions emanating from the cutoff initial data in Section \ref{sec:numerical}. Then, in Section \ref{sec:discussion}, we compare the numerical results to our proposed analytical solution. We finish with concluding remarks in Section \ref{sec:conclusion}.

\section{Construction of the Linear Velocity Solution}\label{sec:linear_vel_sol}

Here we review the linear velocity solutions in \cite{Giron2020} and \cite{Ramsey_linear_velocity} with the concept of linear velocity solutions originating from \cite{Sedov_book}. In doing so, we obtain the homogeneous isentropic implosion solution introduced in \cite{Kidder1974}. When deriving the linear velocity solutions, we will deviate from the derivation done in \cite{Giron2020} and \cite{Ramsey_linear_velocity}, giving a new derivation of this solution.

\par We introduce the \textit{linear velocity ansatz}:
\begin{equation}\label{eq:lin_vel}
    v(r,t) = r \frac{\dot{R}(t)}{R(t)}
\end{equation}
where $R(t)$ is the \textit{scale radius}, to be determined. This velocity is a homogeneous compression as introduced in \cite{Kidder1974}. The scale radius represents the Lagrangian position of a particle \cite{Sedov_book}.
\begin{comment}
We interpret the scale radius as the Lagrangian position of a particle as we see from the following derivation. The position of a particle originating at position $r_0$ is given by 
\begin{equation}
    r(t) = \int_0^t v(s)ds + r(0).
\end{equation}
Substituting the linear velocity ansatz into equation \eqref{eq:lagrangian_position}, we see
\begin{align*}
    r(t) &= \int_0^t r\frac{\dot{R}(s)}{R(s)}ds + r_0 \\
    \frac{dr}{dt} &= r \frac{\dot{R}(t)}{R(t)} \\
    r(t) &= r_0 R(t),
\end{align*}
giving a firm realization to the scale radius $R(t)$.
\end{comment}

%\par Plugging (\ref{eq:lin_vel}) into (\ref{eq:poten_mass}) gives
%\begin{equation}
%    \frac{\partial \rho}{\partial t} + r \frac{\dot{R}}{R} \frac{\partial \rho}{\partial r} + (n+1) \rho \frac{\dot{R}}{R} = 0,
%\end{equation}
%where $n$ is as defined on page 1. 
% The above equation has characteristics
%\begin{equation}\label{eq:char_system_mass}
 %   \frac{dt}{1} = \frac{dr}{r \dot{R}/R} = -\frac{d\rho}{(n+1) (\dot{R}/R) \rho}.
%\end{equation}
%By solving the first and second members of the characteristic system in Equation \eqref{eq:char_system_mass}, we obtain $R(t)\xi = r$ for some $\xi > 0$. By solving the first and third members, we obtain $\rho= R^{-(n+1)}\mathcal{R}(\xi)$. We thus have a solution for the density,
We assume that $\rho$ depends on the similarity variable $\xi = r/ R(t)$.  Conservation of mass (in all of space if $\rho$ is globally integrable or locally inside a ball of radius $MR(t)$ otherwise) implies that $\rho$ has the general form
\begin{equation}\label{eq:density_ss}
    \rho(\xi,t) = \frac{\mathcal{R}(\xi)}{R^{n+1}(t)}
\end{equation}
where $\mathcal{R}(\xi)$ is an arbitrary function of $\xi$. Since we consider the isentropic case, the pressure is $C\rho^\gamma$. We substitute the ansatzes \eqref{eq:lin_vel} and \eqref{eq:density_ss} into \eqref{eq:isen_mom} to obtain
\begin{equation}\label{eq:seperation_of_variables}
    R(t)^{1 - \eta} \ddot{R}(t) = -\frac{C\gamma}{\xi} \mathcal{R}'(\xi) \mathcal{R}(\xi)^{\gamma - 2}
\end{equation}
where $\eta = (n+1)(1 - \gamma)$. By separation of variables, both sides of equation \eqref{eq:seperation_of_variables} evaluate to constants. The left hand side of equation \eqref{eq:seperation_of_variables} implies
\begin{equation}\label{eq:ddotR}
    \frac{d}{dt} \left( \ddot{R}R^{1 - \eta} \right) = 0.
\end{equation}

\par  We select the initial conditions $R(0) = R_0$, $\dot{R}(0) = 0$, and $\ddot{R}(0) = \ddot{R}_0$ and justify them as follows. The first initial condition $R(0) = R_0$ simply states Lagrangian particles all must start at some position. The third initial condition $\ddot{R}(t) = \ddot{R}_0$ indicates whether the particles initially accelerate or decelerate. The second initial condition $\dot{R}(0) = 0$ determines whether the particles are at rest at time $t = 0$.

\par We are interested in the case when $\ddot{R_0}  < 0$, corresponding to the imploding case. Without loss of generality, we pick $R_0 = 1$ and $\ddot{R}_0 = -1$. With the selection of initial conditions for $R(t)$, both sides of equation \eqref{eq:seperation_of_variables} equate to $-1$. The right hand side of equation \eqref{eq:seperation_of_variables} simplifies to
\begin{equation}
    \xi = C \gamma \RR^{\gamma - 2} \frac{d\RR}{d\xi}
\end{equation}
which has solution
\begin{equation}\label{eq:RR}
    \RR(\xi) = \left( B\frac{\gamma - 1}{\gamma} \left( \frac{1}{2}\xi^2 + \rho_0 \right) \right)^{\frac{1}{\gamma - 1}}
\end{equation}
where $B = \frac{1}{C}$. With both $R(t)$ and $\RR(\xi)$ defined, the density and pressure are
\begin{equation}\label{eq:density_RR}
    \rho(r,t) = \frac{\left( B\frac{\gamma - 1}{\gamma} \left( \frac{1}{2}\frac{r^2}{R(t)^2} + \rho_0 \right) \right)^{\frac{1}{\gamma - 1}}}{R(t)^{n+1}}
\end{equation}
and
\begin{equation}\label{eq:pressure_RR}
    P(r,t) = C \frac{\left( B\frac{\gamma - 1}{\gamma} \left( \frac{1}{2}\frac{r^2}{R(t)^2} + \rho_0 \right) \right)^{\frac{\gamma}{\gamma - 1}}}{R(t)^{\gamma(n+1)}}
\end{equation}
respectively. This is a generalization to arbitrary heat capacity ratio and dimension of the isentropic implosion solutions described in \cite{Giron2020} and \cite{Ramsey_linear_velocity}.

\par In the special case when $\eta = -2$, we can solve equation \eqref{eq:ddotR} analytically to obtain 
\begin{equation}\label{eq:R(t)}
    R(t) = \left( 1 - t^2 \right)^{\frac{1}{2}}
\end{equation}
leading to the velocity profile 
\begin{equation}\label{eq:velocity_linear}
    v(r,t) = -\frac{rt}{1-t^2}.
\end{equation}
Note the case of $\eta = -2$ has direct physical application as it represents the case when $n = 2$ and $\gamma = \frac{5}{3}$: a monatomic gas in three dimensions. The case of $\eta = -2$ also corresponds to Kidder's solution. Once $R(t)$ is known, the Lagrangian position of particles $X(t,x_0)$, where $x_0$ represents the initial location of particles is, is given by
\begin{equation}\label{eq:lagrangian}
    X(t,x_0) = x_0\sqrt{1 - t^2}.
\end{equation}
So, all particles arrive at the origin by time $t = 1$.

\par We can understand how the implosion of solutions \eqref{eq:density_RR} and \eqref{eq:pressure_RR} occur through the potential when $R(t)$ is given by \eqref{eq:R(t)}. The potential $\phi(r,t)$ in this case is
\begin{equation}\label{eq:potential_RR}
    \phi(r,t) = -\frac{r^2t}{2(1-t^2)}.
\end{equation}
When $t > 0$, the potential at the origin remains zero while at all other points in the domain, the potential decreases over time. This means the origin is the local (and also global) maximum of the potential. In contrast, the origin is the local (and also global) minimum of the density. From equation \eqref{eq:poten_mom}, the density will tend towards local maxima of the potential, which is located at the origin. However, the potential decreases in value everywhere in the domain except for the origin, which means the origin will always remain the global maximum of the potential. So, the flow will always be directed towards the origin up until implosion time. We emphasize how the origin remains attractive for all time because it continues to be the global minimum in density. This is in contrast to the scenario described for attractive local maxima where the pressure term discourages unbounded accumulation.

\par Regardless of $n$ and $\gamma$, these solutions have the property that $\rho(r,t)$ goes to $\infty$ as $r \to \infty$. This feature is physically undesirable. Although the initial velocity is zero, the initial density is radially unbounded. Phrased differently, the total initial energy of these solutions is radially unbounded because the initial potential energy is radially unbounded. This motivates our investigation of ``cutoff" solutions - solutions emanating from the Kidder solution initial data on some finite spatial interval $[0,x^*]$ and is continuously connected to a constant outside the interval.

\par For the remainder of this paper, we will consider the imploding solution in the specific case of $n = 0$ (one dimension) and $\gamma = 3$ in order to ensure $\eta = -2$, leading to the convenient formula for $R(t)$ prescribed in equation \eqref{eq:R(t)}. We also set the constants $\rho_0,B,C$ described previously as $\rho_0 = 0$ and $B = C = 1$ respectively. This gives initial data
\begin{equation}\label{eq:implosion_initial}
    \boxed{\rho(x,0) = \frac{\sqrt{3}}{3}x, \quad v(x,0)= 0}
\end{equation}
leading to the imploding solution
\begin{equation}\label{eq:implosion_solution}
    \boxed{\rho(x,t) = \frac{\sqrt{3}}{3} \frac{x}{1 - t^2}, \quad v(x,t) = -\frac{xt}{1 - t^2}.}
\end{equation}
As mentioned before, the initial density is radially unbounded. The implosion time is $t^* = 1$, and the density and velocity approach infinity everywhere except the origin as $t \to t^*$. We refer to this problem and corresponding solution as the \textit{Kidder} problem and solution.

\section{1D Cutoff Implosion}\label{sec:1D_cutoff}
In this section, we will construct the solution to the imploding initial data with a constant cutoff in one dimension described in the previous section and discuss some pertinent features of the solution.

\subsection{Introduction of Cutoff Initial Data}
\par In one dimension $\rho$ is symmetric and $v$ is anti-symmetric about the origin $x=0$. As previously mentioned, we are interested in seeing whether the imploding behavior persists even when the initial data in $\rho$ is cut off in the far field. Continuously cutting off the initial data of \eqref{eq:implosion_solution} with a constant yields 
\begin{equation}\label{eq:implosion_cutoff_initial}
\rho(x,0) = \begin{cases}
    \frac{\sqrt{3}}{3}x &\text{ for } |x| \leq x^* \\
    \frac{\sqrt{3}}{3}x^* &\text{ for } |x| \geq x^*    
\end{cases},
\quad v(x,0)= 0
\end{equation}
where $x^*$ is a positive real number. We call this problem and the corresponding solution the \textit{Cutoff Implosion} problem and solution. This initial density is equal to the Kidder initial density for $|x|\leq x^*$ and is a constant outside of this interval. As the solution evolves, new behavior emanates from the cutoff point $x^*$. 
We analyze this behavior in the following section.

\subsection{Hyperbolic Conservation Law Theory}
We first review the necessary hyperbolic conservation law theory as described in \cite{LaxMonograph}, \cite{LevequeGreenBook}, and \cite{EvansPDE} to construct our cutoff implosion solution. We rewrite the one dimensional isentropic Euler equations (\eqref{eq:isen_mass} and \eqref{eq:isen_mom}) in conservation form as
\begin{equation}\label{eq:potential_flow_conservation}
        \rho_t + (\rho v)_x = 0 , \quad
        v_t + \frac{1}{2}(3\rho^2 + v^2)_x = 0
\end{equation}
defining the fluxes as 
\begin{equation}\label{eq:fluxes}
    f(\rho,v) = v\rho, \quad g(\rho,v) = \frac{1}{2} \left( 3\rho^2 + v^2 \right).
\end{equation} We can also write the potential flow equations in matrix form:
\begin{equation}\label{eq:potential_flow_matrix}
    \begin{bmatrix}
        \rho_t \\ v_t
    \end{bmatrix} + \begin{bmatrix}
        v & \rho \\ 3\rho & v
    \end{bmatrix}\begin{bmatrix}
        \rho_x \\ v_x
    \end{bmatrix} = 0.
\end{equation}
Let $\vec{u} = \begin{bmatrix}
    \rho & v
\end{bmatrix}^T$ and $A$ be the coefficient matrix in front of $\begin{bmatrix}
    \rho_x & v_x
\end{bmatrix}^T$. The eigenvalues and corresponding normalized eigenvectors of $A$ are 
\begin{equation}\label{eq:eigenvalues}
    \lambda = v \pm \sqrt{3}\rho, \quad \vec{r}_{\pm} = \frac{1}{2}\begin{bmatrix}
    1 & \pm \sqrt{3}
\end{bmatrix}^T.
\end{equation}
The system is strictly hyperbolic as long as $\rho \neq 0$. Note $\sqrt{3}\rho$ is the sound speed from equation \eqref{eq:sos_potential}. 

\par We verify that the system is genuinely nonlinear by determining whether 
\begin{equation}\label{eq:gen_nonlinear}
    \nabla \lambda_{\pm} \cdot \vec{r}_{\pm} \neq 0,
\end{equation}
which for the system in equation \eqref{eq:potential_flow_matrix}, evaluates to $\pm 2\sqrt{3}$. So the system is genuinely nonlinear. Knowledge of the eigenvalues and the system's genuine nonlinearity is key to construct a rarefaction solution, which we do in the following section.

\subsection{Rarefaction Theory}
\begin{figure}[h]
    \centering    
    \includegraphics[width=0.65\linewidth]{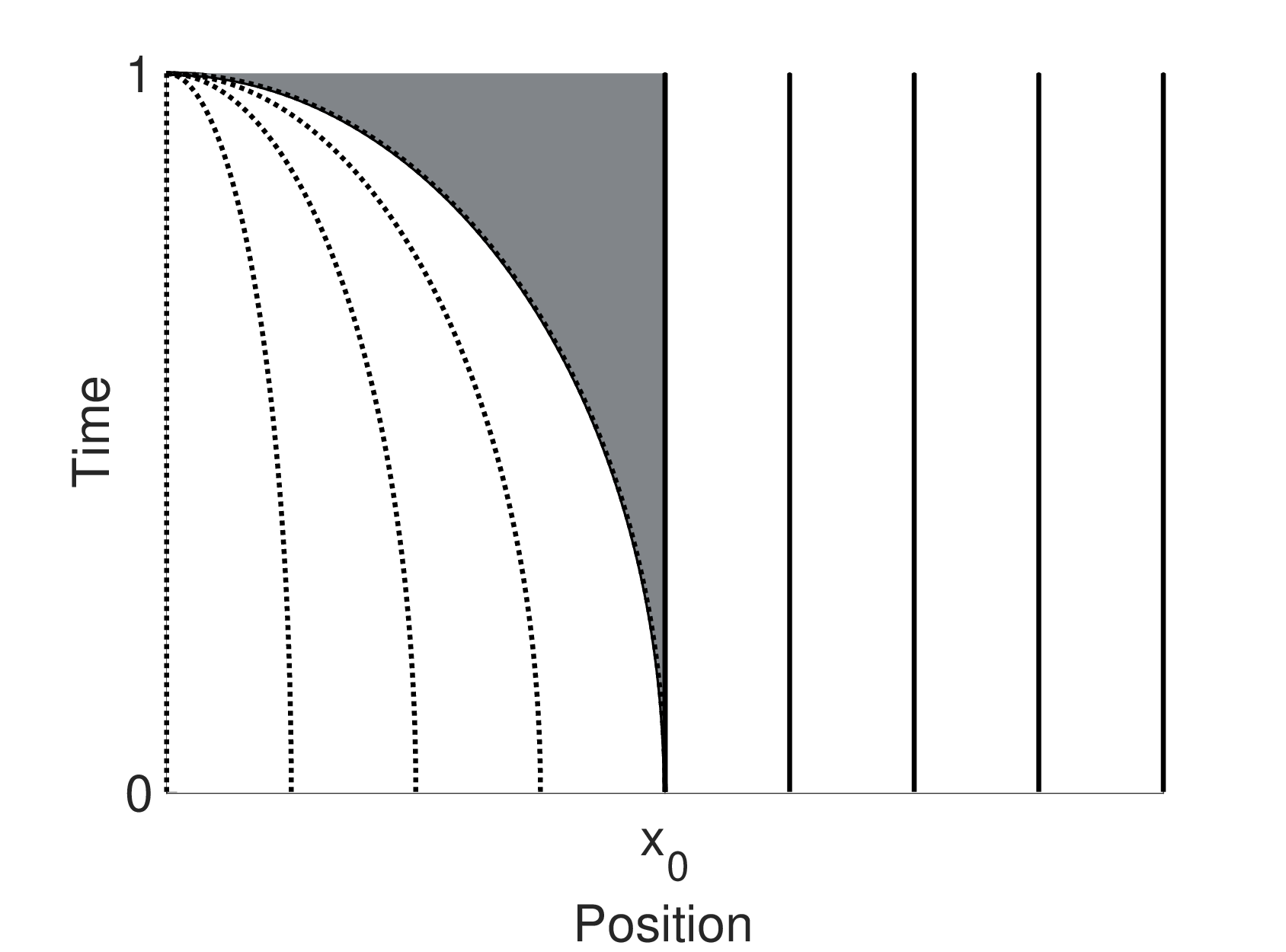}
    \caption{Predicted characteristics stemming from the cutoff implosion initial data. The characteristics left of the cutoff position $x_0$ (dashed lines) are expected to follow the Lagrangian position described in equation \eqref{eq:lagrangian}. The characteristics right of the cutoff position $x_0$ (solid lines) are expected to remain vertical since the initial velocity and density is constant in this region. A rarefaction is expected to emerge from $x_0$ to fill in the empty space of the characteristic diagram, indicated by the region shaded in gray.}
    \label{fig:characterestics}
\end{figure}

We first give a brief heuristic for why we expect a rarefaction to emerge. The conservation of momentum equation \eqref{eq:potential_flow_matrix} can be written as 
\[
v_t + vv_x = -3\rho\rho_x.
\]
which is Burgers' equation with a source term. Although the velocity is initially zero, the source term on the right instantaneously initializes the velocity. But due to the non-differentiability of the cutoff, the source term is initially discontinuous with the source term being negative for all $x < x^*$ and zero for all $x > x^*$. This discontinuity in source term would lead to a discontinuity in velocity as well, causing a rarefaction to emerge from $x^*$. Another way of arriving at this conclusion is by considering the Lagrangian particle positions. For $x < x_0$, the initial data is the same as the Kidder problem initial data. The Lagrangian position of particles for the Kidder solution is given by equation \eqref{eq:lagrangian}. If all particles originating from $x < x_0$ follow this trajectory, they would all arrive at the origin at time $t = 1$. However, since the solution is stationary to the right of $x_0$, the Lagrangian position of those particles remain stationary. The contrasting characteristic behavior on different sides of $x^*$ leads to a gap in the characteristic diagram, indicating a rarefaction should fill in the empty space in the characteristic diagram. A sketch of this characteristic diagram is shown in Figure \ref{fig:characterestics}. We are thus motivated to look for rarefaction solutions to the original problem.

\par To construct a rarefaction solution, suppose the system admits the similarity variable $y = \frac{x}{t}$. The similarity variable may be shifted in space and time arbitrarily due to the time and position shift invariance of systems of hyperbolic conservation laws in 1D planar geometry. Also suppose on the rarefaction region $y_1 \leq y \leq y_2$ the solution is of the form $\vec{w}(y) = \begin{bmatrix}
    \rho(\frac{x}{t+a}) \quad v(\frac{x}{t+a})
\end{bmatrix}^T$ where $y_1, y_2$ is to be determined. The similarity solution $\vec{w}$ must satisfy the ODE
\begin{equation}\label{eq:w_ODE}
    \vec{w}'(y) = \frac{r_{\pm}(y)}{\nabla \lambda_{\pm}(\vec{w}(y)) \cdot r_{\pm}(y)}.
\end{equation}
The denominator of \eqref{eq:w_ODE} will always be nonzero as long as genuine nonlinearity, condition \eqref{eq:gen_nonlinear}, is satisfied.

\subsection{Construction of the Cutoff Implosion Solution}
We solve the ODE system in \eqref{eq:w_ODE} with our eigenvalues and eigenvectors \eqref{eq:eigenvalues} to obtain
\begin{align*}
    \begin{bmatrix}
        \rho'(y) \\ v'(y)
    \end{bmatrix} &= \frac{1}{2} \begin{bmatrix}
    \pm \frac{\sqrt{3}}{3} \\ 1 
    \end{bmatrix}, \\
\end{align*}
which have the solutions
\begin{equation}\label{eq:rare_1}
    \rho(y) = \pm \frac{1}{2} \frac{\sqrt{3}}{3} y + \rho_0, \quad v(y) = \frac{1}{2} y + v_0.
\end{equation}
with the sign dependent on whether we consider $\lambda_+$ or $\lambda_-$.
\par We now address admissibility. First, since we have genuine nonlinearity, $\lambda_{\pm}$ is monotone increasing along the integral curve. This means that two states $u_l$ and $u_r$ can always be connected by a rarefaction wave given
\begin{equation}\label{eq:rare_admit}
    \lambda_{\pm}(u_L) < \lambda_{\pm}(u_R).
\end{equation}

%Now we have all the necessary ingredients to construct the rarefaction solution. 
\par From the previous discussion, we expect to see a rarefaction solution emerging from the cutoff position $x^*$. Arguing via finite speed of propagation for the strictly hyperbolic system, we posit for all times $t < \tilde{t}$ where $\tilde{t}$ is yet to be determined and may be infinite, there exists an interval $[0,\tilde{x}(t))$ for positive real valued function $\tilde{x}(t)$ such that the solution is the Kidder solution on this interval.

\par We wish to determine when the admissibility criteria \eqref{eq:rare_admit} holds with left state $\vec{u}_l = \begin{bmatrix}
    \frac{\sqrt{3}}{3} \frac{x}{1 - t^2} & -\frac{xt}{1-t^2}
\end{bmatrix}^T$ (representing the Kidder solution) and right state $\vec{u}_r = \begin{bmatrix}
    \frac{\sqrt{3}}{3} x^* & 0
\end{bmatrix}^T$ (representing the constant state). By checking the values of $(x,t)$ such that $\lambda_+(\vec{u}_l) < \lambda_+(\vec{u}_r)$ holds with the choices of left and right state, we obtain the inequality
\begin{align*}
    \frac{x}{t+1} &< x^*.
\end{align*}
The 2-eigenvalue $\lambda_+$ is used since this rarefaction corresponds to the 2-rarefaction. Applying the same methodology with the 1-rarefaction emerging from $-x^*$ yields equivalent results. From the above analysis, the upper bound of the rarefaction region is given by $x^*(1+t)$. More provocatively, this implies the similarity variable is $y = \frac{x}{t+1}$, which is a valid form for the similarity variable given the previous comments about the similarity variable permitting time shifts.

\par Continuing with the construction of the rarefaction analysis, the rarefaction solution must continuously agree with the right state at $y = x^*$. Applying that initial condition into equation \eqref{eq:rare_1} in the original coordinates $(x,t)$ yields
\begin{equation}
 \rho(x,t) = \frac{1}{2} \frac{\sqrt{3}}{3} \frac{x}{t+1} + \frac{\sqrt{3}}{3}, \quad     v(x,t) = \frac{1}{2}\frac{x}{t+1} - \frac{1}{2}
\end{equation}
on the rarefaction region.

\par So far, we have determined the form of the rarefaction and the upper bound on the rarefaction region, but such analysis has not indicated the lower bound of the rarefaction region. To obtain this, we argue the rarefaction solution must connect continuously to the Kidder solution at the previously defined position $\tilde{x}(t)$ for all times the implosion solution exists. Thus, we search for $\tilde{x}(t)$ such that the density and velocity match in the Kidder solution and in the rarefaction solution:
\[
\frac{\sqrt{3}}{3} \frac{\tilde{x}(t)}{1-t^2} = \frac{1}{2} \frac{\sqrt{3}}{3} \left( \frac{\tilde{x}(t)}{1+t} + x^* \right) \quad \text{and} \quad -\frac{\tilde{x}(t)t}{1 - t^2} = \frac{1}{2}\frac{\tilde{x}(t)}{1+t} - \frac{1}{2}x^*
\]
to find $\tilde{x}(t) = x^*(t - 1)$. We now have completed our description of the rarefaction and can explicitly state the solution of the cutoff implosion problem. The density is given by 
\begin{equation}\label{eq:density_cutoff_implosion}
\boxed{
    \rho(x,t) = \begin{cases}
        \frac{\sqrt{3}}{3}\frac{|x|}{1 - t^2} & \text{ if } |x| \leq x^*(1-t) \\
        \frac{1}{2}\frac{\sqrt{3}}{3}\left( \frac{|x|}{1+t} + x^* \right) & \text{ if } x^*(1-t) \leq |x| \leq x^*(1+t) \\
        \frac{\sqrt{3}}{3}x^* & \text{ if } x^*(1+t) \leq |x|
    \end{cases},
}
\end{equation}
and the velocity is given by 
\begin{equation}\label{eq:velocity_cutoff_implosion}
\boxed{
    v(x,t) = \begin{cases}
        -\frac{xt}{1 - t^2} & \text{ if } |x| \leq x^*(1-t)\\
        \frac{1}{2} \left(\frac{x}{1+t} - x^*\sign{x} \right) & \text{ if } x^*(1-t) \leq |x| \leq x^*(1+t)\\
        0 & \text{ if } x^*(1+t) \leq |x|
    \end{cases}.
    }
\end{equation}

\subsection{Analysis of the Cutoff Implosion Solution}
We now discuss some of the salient features of our constructed cutoff implosion solution. The most pressing question is whether the cutoff implosion solution still implodes. We have determined the lower bound of the rarefaction region is $x^*(1-t)$. The rarefaction reaches the origin by time $t = 1$, which remarkably is also the original implosion time. On the other hand, the region where the Kidder solution exists on vanishes at time $t = 1$, regardless of where the cutoff is.

\par This does not necessarily imply the cutoff implosion solution does not implode. It may be true that as $t \to 1$, the density around the origin still grows without bound. We will show this cannot be the case via asymptotics. Let $t \to 1$ and set $t = 1 - \epsilon$ for small $\epsilon > 0$. We check the value of the density and velocity at the lower bound of the rarefaction region, so at position $x = x^*(1-t) = \epsilon x^*$. For density, we find
\begin{align*}
    \rho(\epsilon x^*, 1 - \epsilon) &= \frac{\sqrt{3}}{3} \frac{x^*}{2 - \epsilon}
\end{align*}
which approaches $\frac{1}{2}\frac{\sqrt{3}}{3}x^*$ as $\epsilon \to 0$, or equivalently $t \to 1$. For velocity, we find
\begin{align*}
    v(\epsilon x^*, 1 - \epsilon) &= x^* \frac{\epsilon - 1}{2 - \epsilon}
\end{align*}
which approaches $\frac{-x^*}{2}$ as $\epsilon \to 0$. From our asymptotic analysis, we see the rarefaction suppresses the Kidder solution and the cutoff imploding solution does not implode. Again, this is true no matter where the original cutoff is, but we can increase the density near the origin as large as desired by increasing the cutoff position $x^*$.

\par We comment more on the rarefaction itself. The rarefaction variable $y = \frac{x}{t + 1}$ is an unusual modification of the standard rarefaction $\frac{x}{t}$. Although unusual, the time shifted rarefaction variable occurs due to having a non constant left state. This setup is similar to the generalized Riemann problem, where both initial states are allowed to be non constant (see \cite{General_Riemann} for more details). In particular, the time shift of $1$ in the rarefaction variable is notable as $t = 1$ is the original imploding time of the Kidder solution. We find it remarkable that the rarefaction variable encodes the original imploding time although the cutoff implosion solution no longer implodes.

\par The rarefaction region can be viewed as the region of influence originating from the cutoff position. We expect the left and right endpoints of the rarefaction region travels at speed $\lambda_-$ and $\lambda_+$ respectively, and we verify this. Since the solution is a rarefaction and connects continuously to the Kidder solution and constant cutoff to the right, the density and velocity at the left and right endpoints of the rarefaction region are the density and velocity of the Kidder solution and constant cutoff respectively, which we know explicitly. Looking on the positive real axis, the motion of the left endpoint satisfies
\begin{equation}\label{eq:left_endpoint_ode}
    \dot{x}_l = \lambda_-(\rho_l,v_l) = -\frac{x_l}{1 - t}, \quad x_l(0) = x^*
\end{equation}
and the right endpoint satisfies
\begin{equation}\label{eq:right_endpoint_ode}
    \dot{x}_r = \lambda_+(\rho_r,v_r) = x^*, \quad x_r(0) = x^*.
\end{equation}
We can solve ODEs \eqref{eq:left_endpoint_ode} and \eqref{eq:right_endpoint_ode} to obtain the positions for the left endpoint and right endpoint of the rarefaction region as
\[
x_l(t) = x^*(1-t), \quad x_r(t) = x^*(1+t)
\]
which are the rarefaction region endpoints obtained previously. This supports our intuition on how the domain of dependence of the cutoff point appears in the cutoff implosion solution.

\par We end by discussing the speed of sound in the cutoff implosion solution. Noticeably, the solution remains subsonic for all times up to $t = 1$. The speed of sound is obtained by substituting the density in equation \eqref{eq:density_cutoff_implosion} into the speed of sound formula of equation \eqref{eq:sos_potential}:
\begin{equation}\label{eq:sos_cutoff_implosion}
    |c(x,t)| = \begin{cases}
        \frac{x}{1 - t^2} & \text{ if } 0 \leq |x| \leq x^*(1-t) \\
        \frac{1}{2}\left( \frac{x}{1+t} + x^* \right) & \text{ if } x^*(1-t)\leq |x| \leq x^*(1+t) \\
        x^* & \text{ if } x^*(1+t) \leq |x| 
    \end{cases}.
\end{equation}
We compare the magnitude of the velocity to the local speed of sound on every region. Without loss of generality, we consider only the positive real axis. On the Kidder solution region,
\[
\left| -\frac{xt}{1 - t^2} \right| \leq \frac{x}{1 - t^2}
\]
as long as $0 \leq t < 1$, so the imploding solution is subsonic. On the rarefaction region,
\[
 v(x,t) = \frac{1}{2} \left( x^* - \frac{x}{1+t} \right) \leq \frac{1}{2}\left( x^* + \frac{x}{1+t} \right) \leq c(x,t).
\]
Finally, in the constant region the flow is stationary. So, for all time $t < 1$, the flow is subsonic. We postulate whether the property of the flow being subsonic is an indicator that the solution will not implode. This hypothesis is motivated by preliminary numerical investigations for higher dimensions, where the final state of the solution at $t = 1$ is ambiguous and the velocity near the origin becomes supersonic.

\subsection{Uniqueness}
\par We discuss questions about existence and uniqueness in this subsection by appealing to the Weak-Strong Uniqueness results introduced by \cite{Dafermos_entropy} and \cite{Diperna_entropy} and exposited in \cite{Dafermos_book}. The relevant definitions are briefly introduced before stating the result.  

\begin{comment}We wish to be certain that our constructed solution and corresponding methodology is the correct and unique solution stemming from the cutoff implosion initial data. Our worries are assuaged by appealing to the Weak-Strong Uniqueness results introduced by \cite{Dafermos_entropy} and \cite{Diperna_entropy} and exposited in \cite{Dafermos_book}. We briefly introduce the relevant definitions before stating the result. \textcolor{red}{this also might be appendix material}
\end{comment}

\par The concept of entropy-flux pairs for systems of hyperbolic conservation laws is rooted in the thermodynamic entropy of gas dynamics. For a system of $m$ scalar hyperbolic conservation laws 
\begin{equation}\label{eq:conservation_law}
    u_t + F(u)_x = 0
\end{equation}
where $u \in \mathbb{R}^m$, the functions $\Phi,\Psi$ mapping $\mathbb{R}^m$ into $\mathbb{R}$ are an \textit{entropy} and \textit{entropy flux} respectively if $\Phi$ is a convex function with $D\Phi(z) DF(z) = D\Psi(z)$ for $z \in \mathbb{R}^m$. With this satisfied, when $u$ is a solution to equation \eqref{eq:conservation_law}, the entropy and entropy flux satisfy the entropy inequality
\begin{equation}\label{eq:entropy_con_law}
    \Phi(u)_t + \Psi(u)_x \leq 0.
\end{equation}
Equality is achieved when $u$ is a sufficiently smooth solution of equation \eqref{eq:conservation_law}. The entropy stemming from hyperbolic conservation law theory is related to thermodynamic entropy as in many physical systems, the entropy of hyperbolic conservation law theory is taken to be the negative of thermodynamic entropy. The entropy inequality \eqref{eq:entropy_con_law} also restricts which shock solutions are admissible, analogous to the second law of thermodynamics enforcing entropy to be non decreasing. For the isentropic Euler equations, the total energy $E = \frac{1}{2}\rho v^2 + \frac{C}{\gamma - 1}\rho^{\gamma}$ serves as the entropy with entropy flux $Ev + C\rho^{\gamma}v$. Note the total energy is not strictly convex as it loses convexity at vacuum, $\rho = 0$. However, in the specific case of $\gamma = 3$, we employ the entropy $\Phi(\rho,v) = \frac{1}{2}v^2 + \frac{3}{2}\rho^2$ and corresponding entropy flux $\Psi(\rho,v) = \frac{1}{3}v^3 + 3\rho^2 v$ which is strictly convex everywhere.

\par Returning to Weak-Strong Uniqueness, the theorem of \cite{Dafermos_entropy} and \cite{Diperna_entropy} states if \eqref{eq:conservation_law} admits an entropy flux pair with strictly convex entropy, then a locally Lipschitz continuous solution is the unique weak solution of the PDE with corresponding initial condition. By using the previously described entropy pair $(\Phi,\Psi)$ and given the constructed solution in equations \eqref{eq:density_cutoff_implosion} and \eqref{eq:velocity_cutoff_implosion} is locally Lipschitz on all of $\mathbb{R}$ for $t \in [0, 1)$, we conclude our constructed rarefaction solution is the unique weak solution on this time interval.

%%This is a paragraph where we couldn't use relative entropy because of non strict convexity of the entropy
%Although our constructed solution is also locally Lipschitz, we cannot directly apply the result due to the total energy not being strictly convex on the whole domain. In particular, the total energy loses strict convexity at the origin where the density remains zero for all $t < 1$. We work around this issue by applying the result to a neighborhood around the cutoff point. Specifically, for any fixed time $T < 1$, on the domain $[x^* - x^*T, x^* + x^*T]$ the density in the rarefaction solution is bounded from below by a constant depending on $T$ but strictly greater than zero. This avoids any issues with loss of strict convexity of the associated entropy, and we may directly apply the theorem thereby stating the rarefaction solution is the unique weak solution up to $T < 1$ on $[x^* - xT, x^* + xT]$. While this discussion is not a rigorous treatment of uniqueness, it does suggest the constructed solution to the cutoff implosion problem is the unique solution from the cutoff initial data.

\subsection{Continuation Past $t = 1$}
As the cutoff implosion solution does not implode at $t = 1$, we can try to extend the solution past $t = 1$. The previously performed asymptotic analysis describes the behavior of the density and velocity at the origin. At $t = 1$, the density at the origin is continuous but jumps up to $\frac{1}{2}\frac{\sqrt{3}}{3}x^*$ after being zero for all times beforehand. On the other hand, the velocity at the origin is discontinuous at $t = 1$ taking the value of $\frac{x^*}{2}$ from the left and $-\frac{x^*}{2}$ from the right. We can define the initial condition naturally arising from continuing the cutoff implosion solution past $t = 1$ as 
\begin{equation}\label{eq:ic_continue_density}
    \rho(x,1) = \begin{cases}
        \frac{1}{2} \frac{\sqrt{3}}{3} \left( \frac{|x|}{2} + x^* \right) & \text{ if } 0 \leq |x| \leq 2x^* \\
        \frac{\sqrt{3}}{3} x^* & \text{ if } 2x^* \leq |x| \\
    \end{cases}
\end{equation}
and
\begin{equation}\label{eq:ic_continue_velocity}
    v(x,1) = \begin{cases}
        \frac{1}{2}\left(\frac{x}{2} - x^*\sign{x}\right) & \text{ if } 0 < |x| < 2x^* \\
        0 & \text{ if } 2x^* \leq |x|
    \end{cases}.    
\end{equation}
To continue the solution for $t > 1$, we construct a solution to the generalized Riemann problem with initial data given by equations \eqref{eq:ic_continue_density} and \eqref{eq:ic_continue_velocity}. Existence theory showing how the solution of the generalized Riemann problems takes the form of its associated Riemann problem is established in \cite{General_Riemann}. Asymptotic expansions in the similarity variable are considered in \cite{LeFloch_part_1} and constructed explicitly for first order in the case of gas dynamics in \cite{LeFloch_part_2}. The authors of \cite{Ben_artzi_1984} considered generalized Riemann problems with linear initial data to construct higher order numerical schemes, and an extension for arbitrary polynomial initial data is discussed in \cite{Toro2009}. 

\par As seen in the velocity initial condition \eqref{eq:ic_continue_velocity}, the characteristics at the origin at time $t = 1$ collide, forming a shock. Due to the symmetry of the solution, we predict a double shock solution to form. If a shock forms and moves to the right from the origin, it must also form and move to the left from the origin, creating a double shock. The properties of these shocks, such as their position and speed, are reflected as well. To determine the double shock structure and intermediate state emanating from the origin, we treat the initial condition as a Riemann problem locally around the origin at time $t = 1$. 

\par Turning to the theory of Riemann problems, a shock is admissible between two given states $\vec{u_l} = \begin{bmatrix}\rho_l(x,t) & v_l(x,t)\end{bmatrix}^T$ and $\vec{u_r} = \begin{bmatrix}\rho_r(x,t) & v_r(x,t)\end{bmatrix}^T$ if it satisfies the Rankine-Hugoniot condition
\begin{equation}\label{eq:RH}
    \dot{s}(t) = \frac{f(\rho_r,v_r) - f(\rho_l,v_l)}{\rho_r - \rho_l} = \frac{g(\rho_r,v_r) - g(\rho_l,v_l)}{v_r - v_l} 
\end{equation}
with the fluxes $f,g$ as defined in equation \eqref{eq:fluxes}, $s(t)$ denoting the shock position, $\dot{s}(t)$ denoting the shock speed, and the understanding that the left and right states vary in space and time. We note our choice of continuing to use the fluxes in equation \eqref{eq:fluxes} is for continuity and mathematical tractability of this work rather than physical accuracy. In general, we do not expect the Rankine-Hugoniot condition to be satisfied in general between any generic pair of left and right states. So, we search for an intermediate state $u^* = (\rho^*,v^*)$ such that the Rankine-Hugoniot condition is satisfied between $u_l$ and $u^*$, 
\begin{equation}\label{eq:RH_left}
    \frac{f(\rho^*,v^*) - f(\rho_l,v_l)}{\rho^* - \rho_l} = \frac{g(\rho^*,v^*) - g(\rho_l,v_l)}{v^* - v_l}, 
\end{equation}
and $u^*$ and $u_r$,
\begin{equation}\label{eq:RH_right}
    \frac{f(\rho_r,v_r) - f(\rho^*,v^*)}{\rho_r - \rho^*} = \frac{g(\rho_r,v_r) - g(\rho^*,v^*)}{v_r - v^*}. 
\end{equation}
However, not all shock connections satisfying equations \eqref{eq:RH_left} and \eqref{eq:RH_right} are physically admissible, so we also enforce the entropy condition
\begin{equation}\label{eq:entropy_condition}
    \lambda_i(\vec{u}_l(s_i(t),t)) > \dot{s}_i > \lambda_i(\vec{u}_r(s_i(t),t)).
\end{equation}
We may form a local Riemann problem by defining the left and right constant states to be 
$\vec{u_l}^0 = \lim_{x \to 0^-} \begin{bmatrix}
    \rho(x,1) & v(x,1)
\end{bmatrix} = \begin{bmatrix}
    \frac{1}{2}\frac{\sqrt{3}}{3} x^* & \frac{1}{2} x^*
\end{bmatrix}$ and the right constant to be $\vec{u_r}^0 = \lim_{x \to 0^+} \begin{bmatrix}
    \rho(x,1) & v(x,1)
\end{bmatrix} = \begin{bmatrix}
    \frac{1}{2}\frac{\sqrt{3}}{3} x^* & -\frac{1}{2} x^*
\end{bmatrix}$. The solution satisfying the entropy condition \eqref{eq:entropy_condition} has a double shock structure with intermediate state $u^{*,0} = \begin{bmatrix}\rho^{*,0} & v^{*,0}\end{bmatrix}^T = \begin{bmatrix}\frac{\sqrt{3}}{3}x^* & 0\end{bmatrix}^T$. The local Riemann solution acts as the leading order behavior of the solution to the generalized Riemann problem as described in \cite{General_Riemann}, \cite{LeFloch_part_1}, \cite{Ben_artzi_1984}, and \cite{Toro2009} and we refer to it as the \textit{associated Riemann problem}.

\par Now for the generalized Riemann problem, the left and right states are no longer constant, and \textit{a priori} neither is the intermediate state. Let $s(t)$ be the position of the shock on the positive axis with $-s(t)$ being the position of the shock on the negative axis. For the continuation problem, the left state is 
\begin{equation}\label{eq:left_state}
    \vec{u_l}(-s(t),t) = \begin{bmatrix}\frac{1}{2} \frac{\sqrt{3}}{3} \left( \frac{s(t)}{1+t} + x^* \right) &\frac{1}{2}\left(-\frac{s(t)}{1+t} + x^*\right) \end{bmatrix}^T
\end{equation}
and the right state is 
\begin{equation}\label{eq:right_state}
\vec{u_r}(s(t),t) = \begin{bmatrix}\frac{1}{2} \frac{\sqrt{3}}{3} \left( \frac{s(t)}{1+t} + x^* \right) &\frac{1}{2}\left(\frac{s(t)}{1+t} - x^*\right) \end{bmatrix}^T.
\end{equation} 
The left and right states are understood to be the continuation of the rarefaction in equations \eqref{eq:density_cutoff_implosion} and \eqref{eq:velocity_cutoff_implosion}. For brevity, we notate the right density and velocity as $\rho$ and $v$ and the left density and velocity as $\rho$ and $-v$ using the symmetry. Similarly, denote $\rho^* = \rho^*(s) = \rho^*(-s)$ and $v^* = v^*(s) = -v^*(-s)$. Applying the left and right states to equations \eqref{eq:RH_left} and \eqref{eq:RH_right} with the fluxes defined in equation \eqref{eq:fluxes} yield
\begin{equation}\label{eq:RH_left_2}
    \frac{-\rho^* v^* + \rho v}{\rho^* - \rho} = \frac{1}{2}\frac{3(\rho^*)^2 + (v^*)^2 - 3\rho^2 - v^2}{-v^* + v} 
\end{equation}
and
\begin{equation}\label{eq:RH_right_2}
    \frac{\rho v - \rho^*v^*}{\rho - \rho^*} = \frac{1}{2}\frac{3\rho^2 + v^2 - 3(\rho^*)^2 - (v^*)^2}{v - v^*} 
\end{equation}
where $\rho^*$ and $v^*$ are understood to be functions of $s$. Equations \eqref{eq:RH_left_2} and \eqref{eq:RH_right_2} represent the same value with opposite sign because each represents the shock speed to the left and right respectively. Without loss of generality, we solve the Rankine-Hugoniot condition \eqref{eq:RH_right_2} for $v^*$ yielding the polynomial
\begin{equation}\label{eq:v_star_poly}
    0 = (\rho + \rho^**)(v^*)^2 - 2(\rho + \rho^*)v^* - 3(\rho^*)^3 + 3\rho(\rho^*)^2 + (v^2 + 3\rho^2)\rho^* - 3\rho^3 + v^2\rho
\end{equation}
which has roots
\begin{equation}\label{eq:v_star_roots}
    v^* = v \pm \sqrt{3}(\rho^* - \rho)
\end{equation}
using the fact $\rho = \rho^*$ is a root in $\rho$ of the constant term in $v^*$. To determine which of \eqref{eq:v_star_roots} is permissible, the entropy condition \eqref{eq:entropy_condition} is employed to check $\lambda_+(u^*) > \lambda_+(\vec{u}_r)$, finding $v^* = v + \sqrt{3}(\rho^* - \rho)$ is the entropy satisfying intermediate state.

\par In general, the intermediate state on the region between the shocks is not constant, although the constant intermediate state for the associated Riemann problem does satisfy constraint \eqref{eq:v_star_roots}. We will now argue the constant state $u^{*,0}$ is also the intermediate state for the generalized Riemann problem. First, we consider a power series expansion in time around $x = 0, t = 1$ for the solution $\vec{u}(x,t)$ as done in \cite{Toro2009} which takes the form
\begin{equation}\label{eq:u_power_series}
    \vec{u}_{LR}(t) = u^{*,0} + \sum_{k = 1}^K \partial^k_t \vec{u}(0,1^+)\frac{(t-1)^k}{k!}
\end{equation}
for some time $t \geq 1$. 

\par To compute $\partial^k_t \vec{u}(0,1^+)$, \cite{Ben_artzi_1984} and \cite{Toro2009} propose studying the evolution of the spatial derivatives and transforming them into $\partial^k_t \vec{u}(0,1^+)$. In our case, we can directly use the explicit form of the solution in equations \eqref{eq:density_cutoff_implosion} and \eqref{eq:velocity_cutoff_implosion} to obtain the time derivatives, computing $\partial^k_t \vec{u}(0,1^+)$ to be zero. The power series representation is subsequently
\[
\vec{u}_{LR}(t) = u^{*,0},
\]
implying at the origin the solution is the constant intermediate state of the associated Riemann problem for some time interval. Note following the procedure of \cite{Ben_artzi_1984} and \cite{Toro2009} yields the same result. By evaluating the PDEs \eqref{eq:potential_flow_conservation} at $x = 0$ and using the facts $\rho_t(0,t) = v_t(0,t) = 0$ and $v^{*,0} = 0$ on this time interval, it can be shown $\rho_x(0,t) = v_x(0,t) = 0$. Performing the same procedure on PDEs \eqref{eq:potential_flow_conservation} after taking further time derivatives shows $\partial_x^k \rho (0,t) = \partial_x^k v (0,t) = 0$ as well. If $\rho$ and $v$ are assumed to be analytic in space at $x = 0$ on this time interval, then $\rho$ and $v$ are constants on the intermediate region. For this reason, we take the intermediate state for the generalized Riemann problem to be $u^* = u^{*,0} = \begin{bmatrix}\frac{\sqrt{3}}{3}x^* & 0\end{bmatrix}^T$.

\par Thus, the solution takes the form \begin{equation}\label{eq:ic_continue_density_int}
    \rho(x,t) = \begin{cases}
        \frac{\sqrt{3}}{3}x^* & \text{ if } -s(t) \leq |x| \leq s(t) \\
        \frac{1}{2} \frac{\sqrt{3}}{3} \left( \frac{|x|}{1+t} + x^* \right) & \text{ if } s(t) \leq |x| \leq x^*(1+t) \\
        \frac{\sqrt{3}}{3} x^* & \text{ if } x^*(1+t) \leq |x| \\
    \end{cases}
\end{equation}
and
\begin{equation}\label{eq:ic_continue_velocity_int}
    v(x,t) = \begin{cases}
        0 & \text{ if } -s(t) \leq |x| \leq s(t) \\
        \frac{1}{2}(\frac{x}{1+t} - x^*\sign{x^*}) & \text{ if } s(t) \leq |x| \leq x^*(1+t) \\
        0 & \text{ if } x^*(1+t) \leq |x| \\
    \end{cases} 
\end{equation}
All that is left is to determine the position and speed of the shock $s(t)$ emanating from the origin. We apply the Rankine-Hugoniot condition for the right moving shock with left state $\vec{u_l} = \begin{bmatrix}\frac{\sqrt{3}}{3}x^* & 0 \end{bmatrix}^T$ and $\vec{u_r} = \begin{bmatrix} \frac{1}{2} \frac{\sqrt{3}}{3} \left( \frac{s(t)}{1+t} + x^* \right),\frac{1}{2}\left( \frac{s(t)}{1+t} - x^*\right) \end{bmatrix}^T$ to obtain
\begin{equation}\label{eq:shock_speed}
    \dot{s}(t) = \frac{1}{2}\left( \frac{s}{1+t} + x^* \right)
\end{equation}
which is an ODE in $s(t)$. The ODE of equation \eqref{eq:shock_speed} is equipped with the initial condition $s(1) = 0$, signifying the shock originates at the origin at time $t = 1$ and has solution
\begin{equation}\label{eq:shock}
    s(t) = x^* \left((1+t) - (2(1+t))^{\frac{1}{2}} \right).
\end{equation}
The shock cannot catch up to the rarefaction as the right endpoint position of the rarefaction region is given by $x^*(1+t)$. Thus, the solution to the cutoff implosion problem past $t = 1$ is given by equations \eqref{eq:ic_continue_density_int} and \eqref{eq:ic_continue_velocity_int} with shock position given by equation \eqref{eq:shock}.

\section{Numerical Results}\label{sec:numerical}
We complement our analysis with numerical simulations of the potential flow equations in equations \eqref{eq:potential_flow_conservation} with the cutoff implosion initial data in equation \eqref{eq:implosion_cutoff_initial}. Our goal is to verify the analytical solution by comparing it to a numerically computed solution. We also would like to determine whether various predicted features - such as the similarity variable $y$ - are also present in the numerical solution. 

\par To perform the numerical simulation, we use a standard Lagrangian leapfrog scheme as introduced in \cite{VonNeumannArtificialViscosity} and further developed by \cite{NohArtificialViscosity}. In this Lagrangian scheme, the spatial domain is interspersed with nodes, where the kinematic variables such as velocity and position are defined. In between the nodes are the fluid elements, where variables such as the density and pressure are defined. Lagrangian schemes are characterized by how they allow the nodes themselves to move in time. The most prominent feature of the scheme presented in \cite{VonNeumannArtificialViscosity} and \cite{NohArtificialViscosity} is the introduction of artificial viscosity in the pressure terms to better resolve shocks. However, since we do not expect or see shocks in our numerical simulations, we set artificial viscosity terms to zero in our calculations.

\par For the numerical simulations, we use a numerical domain of $[0,5]$ with $1501$ nodes and $1500$ zones. Initially, the nodes are equally spaced, but their positions change over time via the dynamics of the Lagrangian numerical scheme. We use dynamic time stepping to ensure our time steps always meet a Courant–Friedrichs–Lewy (CFL) condition up to ending time $T = 0.999$ with Courant number $C = 0.1$. We set the leftmost node boundary condition as $v = 0$ at $x = 0$ and the rightmost node boundary condition as $v = 0$ at $x = 5$. This is justified as the velocity of the Kidder solution remains $0$ at $x = 0$ up until the original implosion time. For the right boundary condition, as long as the density remains constant at $x = 5$, the velocity will also remain $0$ as well. We select the cutoff position in our initial data to be $x^* = 2.5$.

\par
\begin{figure}[h]
     \centering
     \begin{subfigure}[h]{0.49\textwidth}
         \centering
         \includegraphics[width=\textwidth]{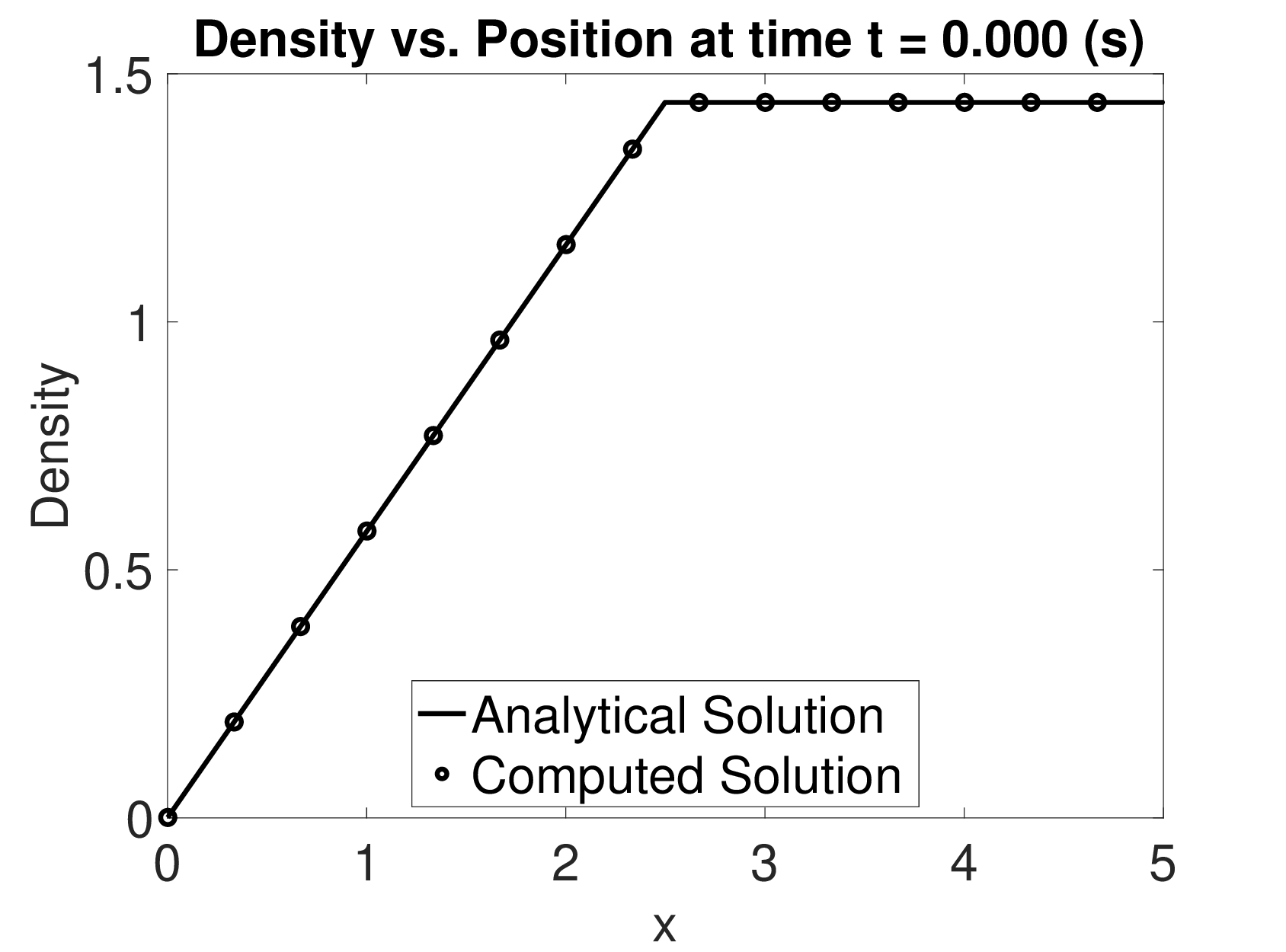}
         \caption{Density}
         \label{subfig:density_ic}
     \end{subfigure}
     \hfill
     \begin{subfigure}[h]{0.49\textwidth}
         \centering
         \includegraphics[width=\textwidth]{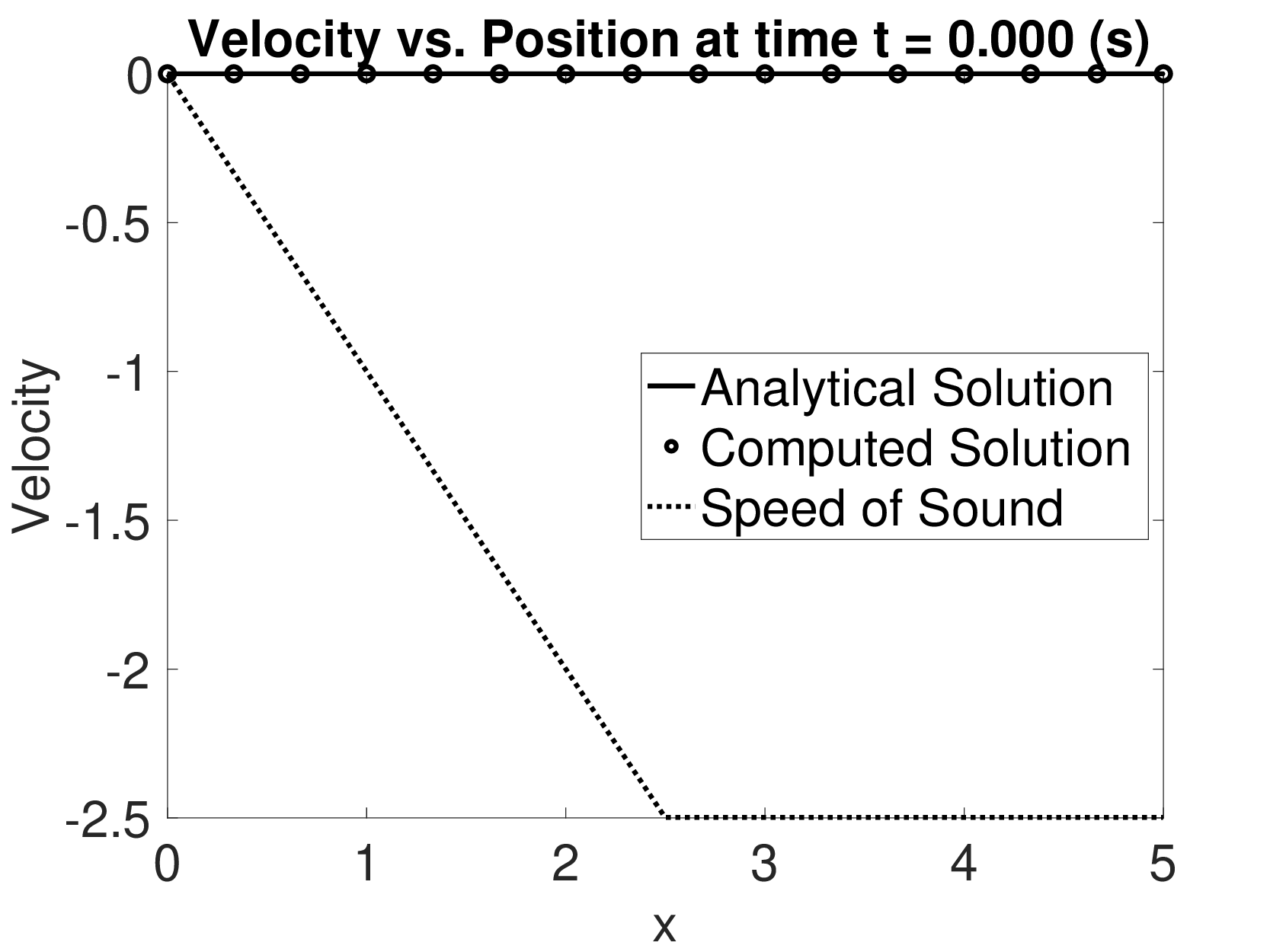}
         \caption{Velocity}
         \label{subfig:velocity_ic}
     \end{subfigure}
        \caption{Initial data for the cutoff implosion problem.}
        \label{fig:initial}
\end{figure}

\begin{figure}[h]
     \centering
     \begin{subfigure}[h]{0.49\textwidth}
         \centering
         \includegraphics[width=\textwidth]{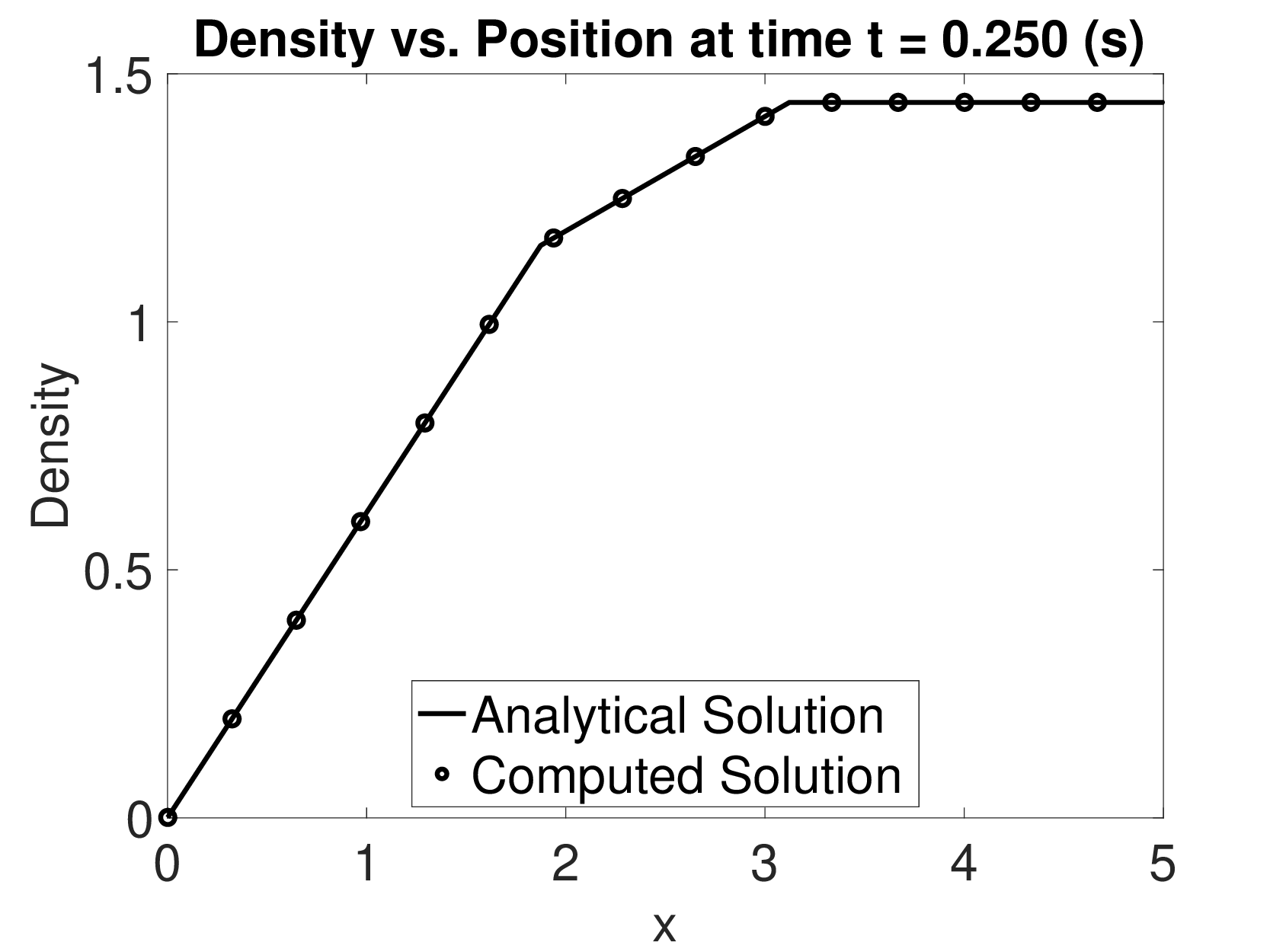}
         \caption{$t \approx 0.25$}
         \label{subfig:density_0.25}
     \end{subfigure}
     \hfill
     \begin{subfigure}[h]{0.49\textwidth}
         \centering
         \includegraphics[width=\textwidth]{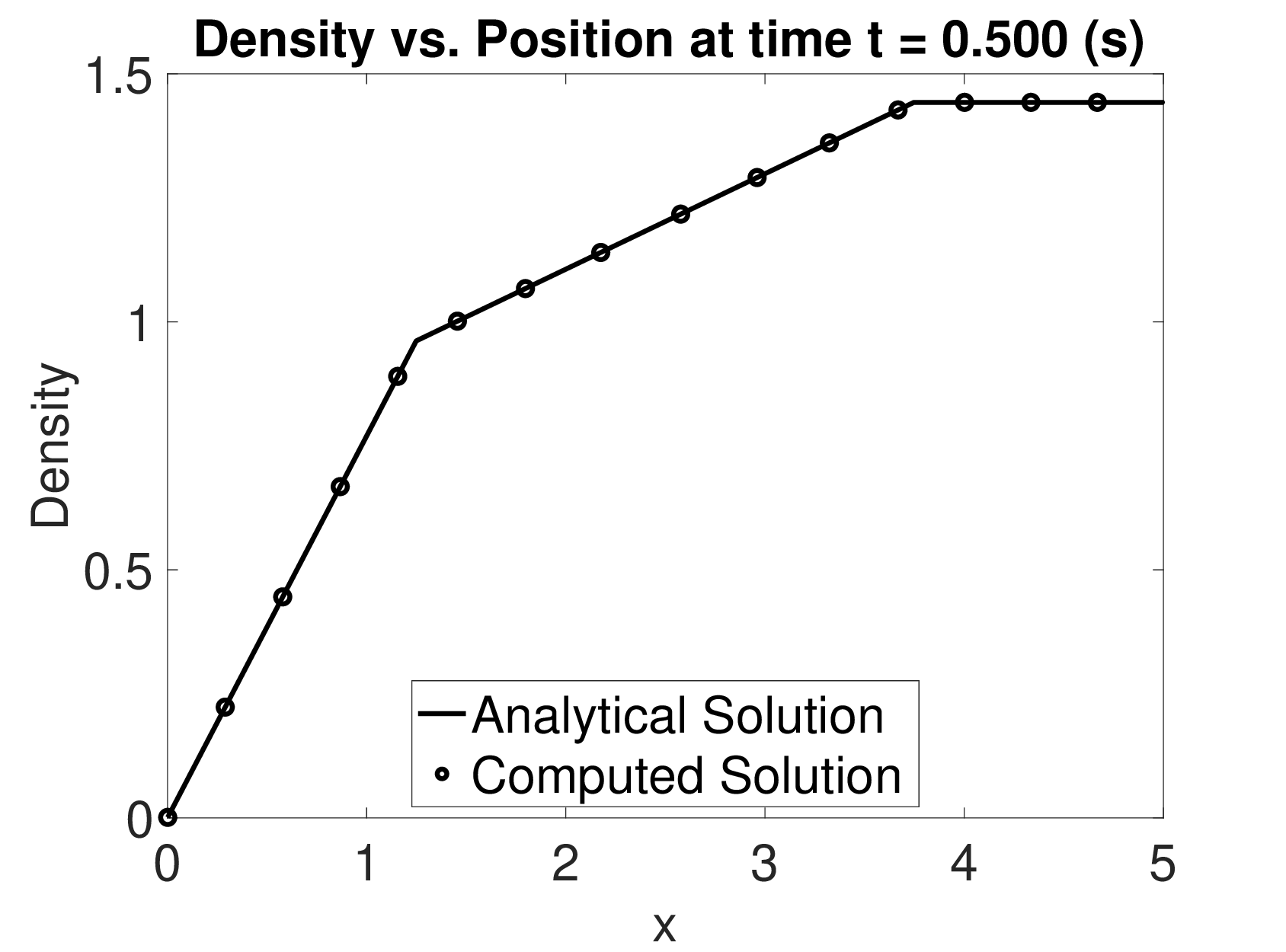}
         \caption{$t \approx 0.5$}
         \label{subfig:density_0.5}
    \end{subfigure}
    \vskip\baselineskip
    \begin{subfigure}[h]{0.49\textwidth}
         \centering
         \includegraphics[width=\textwidth]{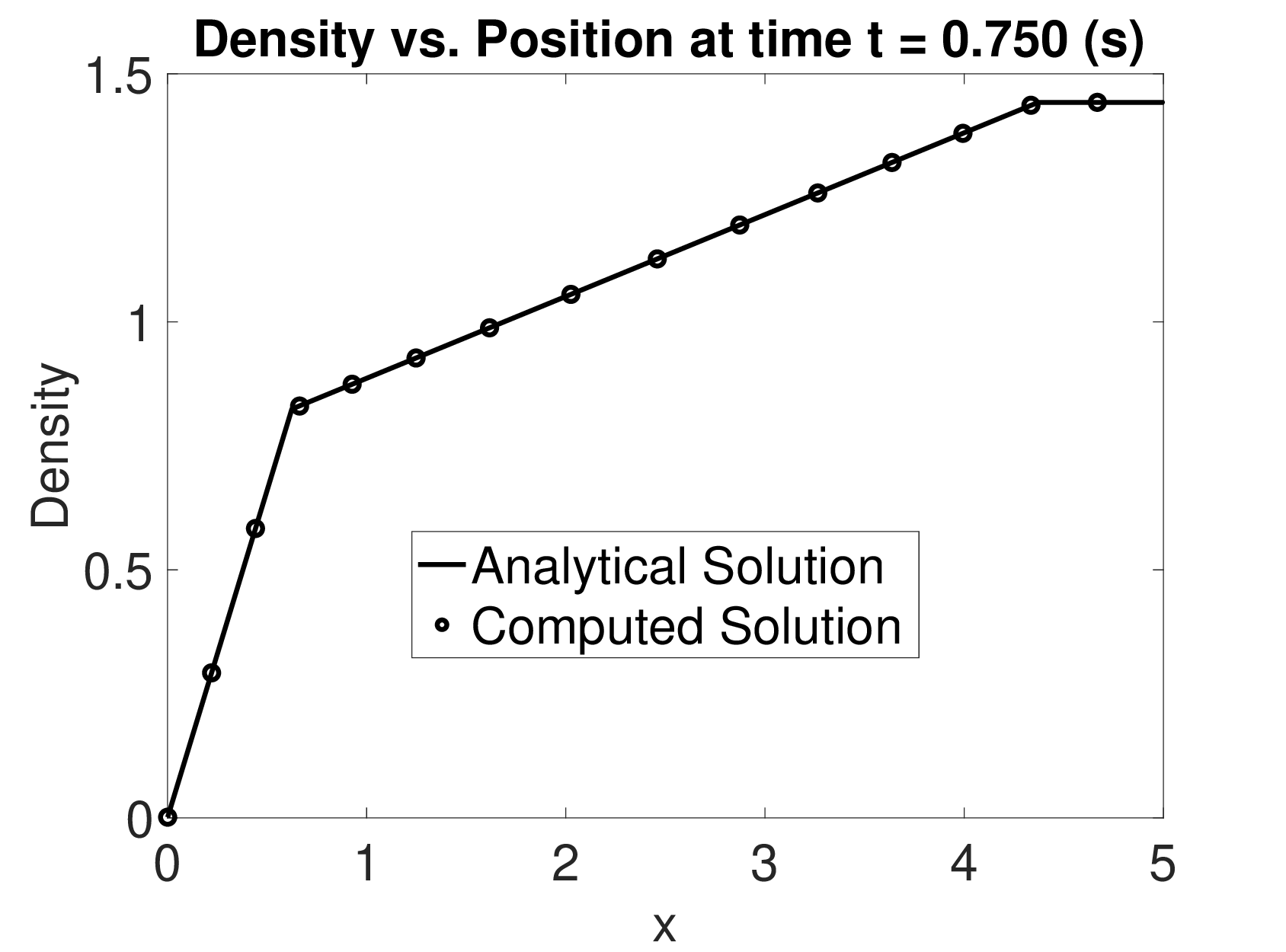}
         \caption{$t \approx 0.75$}
         \label{subfig:density_0.75}
     \end{subfigure}
     \hfill
     \begin{subfigure}[h]{0.49\textwidth}
         \centering
         \includegraphics[width=\textwidth]{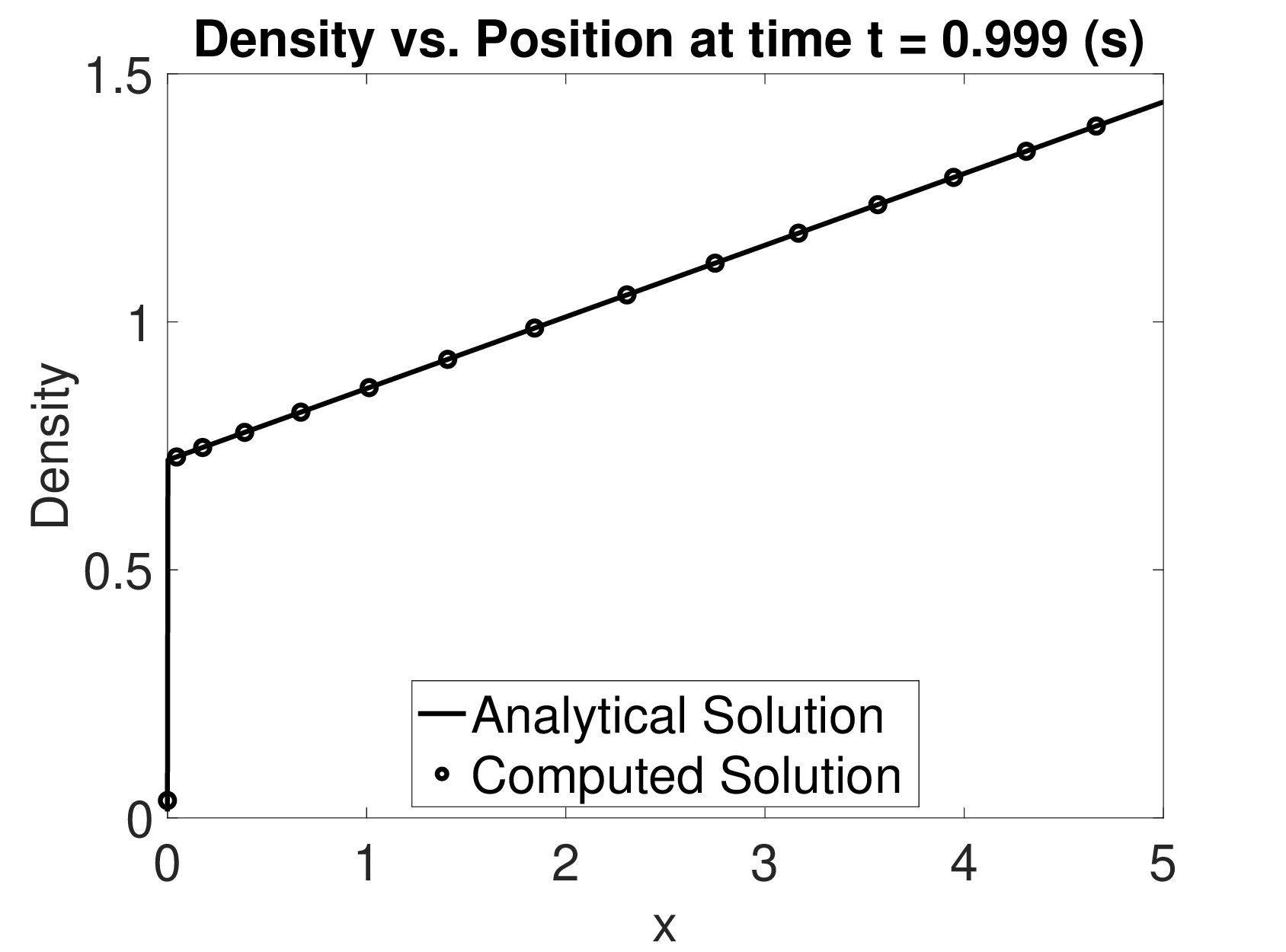}
         \caption{$t \approx 1$}
         \label{subfig:density_1}
     \end{subfigure}
        \caption{Analytical solution and selected representatives of the computed solution for density at approximate times $t = 0.25,0.5,0.75$ and $t = 1$.}
        \label{fig:density}
\end{figure}

\begin{figure}[h]
     \centering
     \begin{subfigure}[h]{0.49\textwidth}
         \centering
         \includegraphics[width=\textwidth]{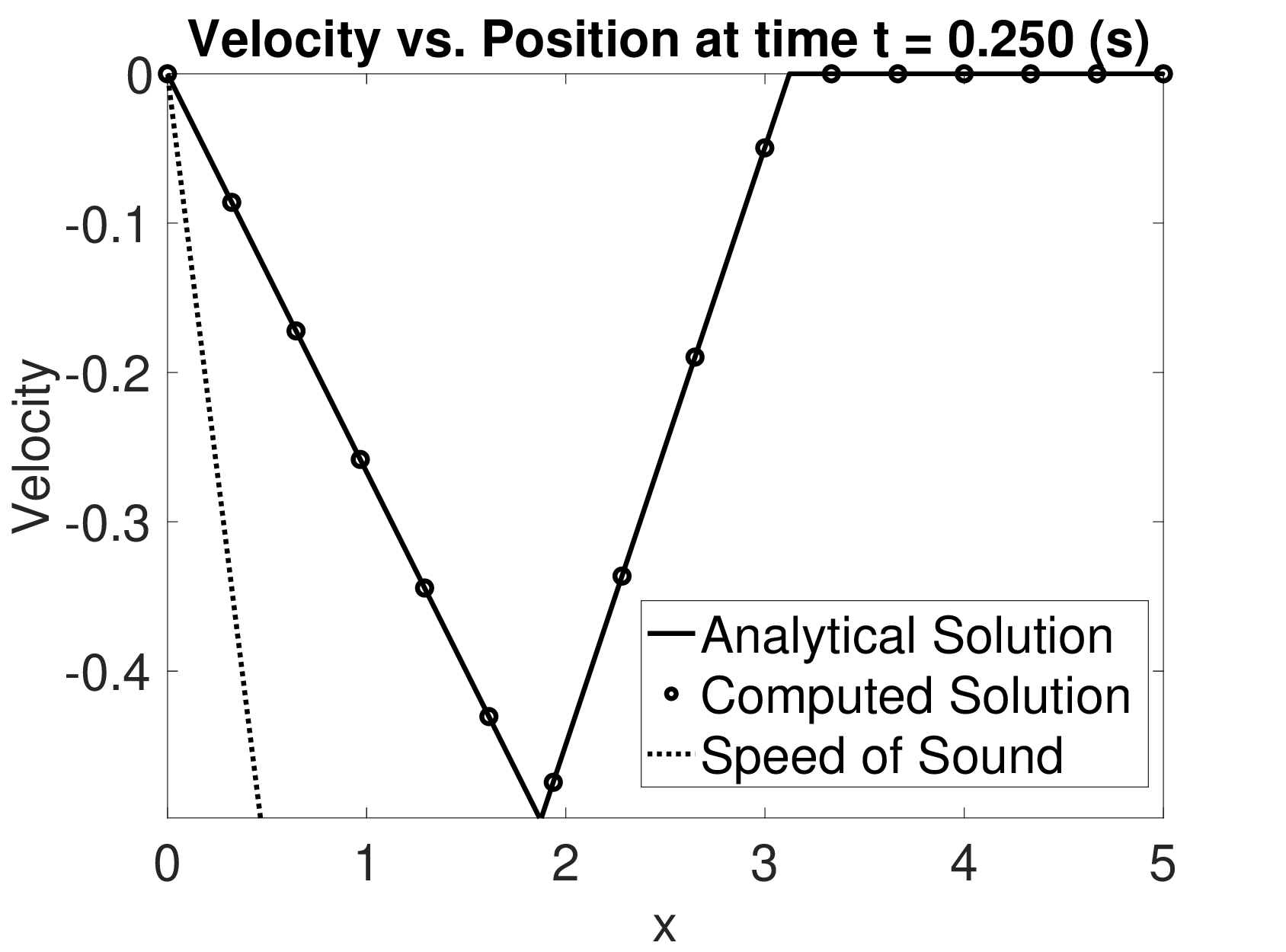}
         \caption{$t \approx 0.25$}
         \label{subfig:velocity_0.25}
     \end{subfigure}
     \hfill
     \begin{subfigure}[h]{0.49\textwidth}
         \centering
         \includegraphics[width=\textwidth]{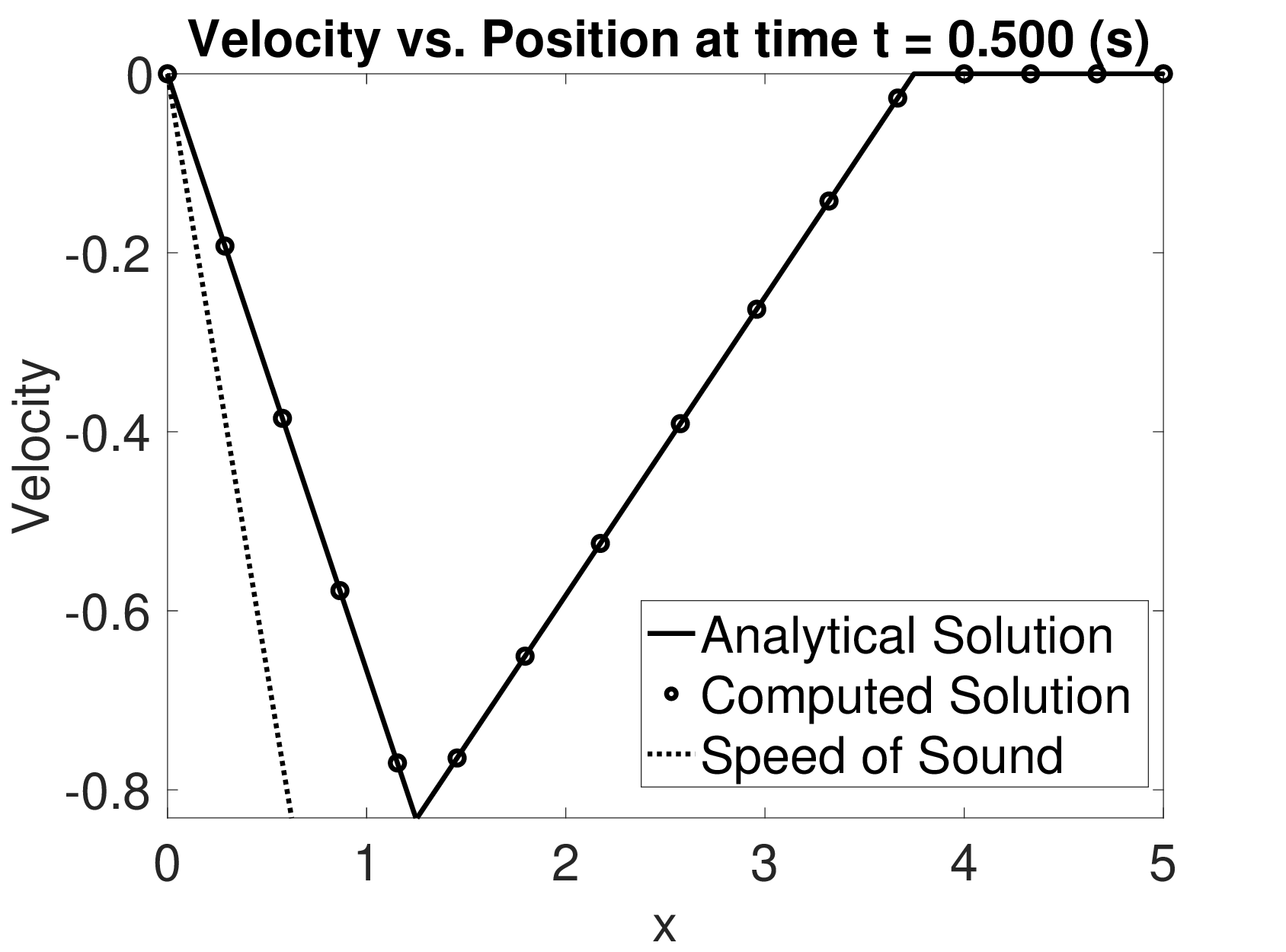}
         \caption{$t \approx 0.5$}
         \label{subfig:velocity_0.5}
    \end{subfigure}
    \vskip\baselineskip
    \begin{subfigure}[h]{0.49\textwidth}
         \centering
         \includegraphics[width=\textwidth]{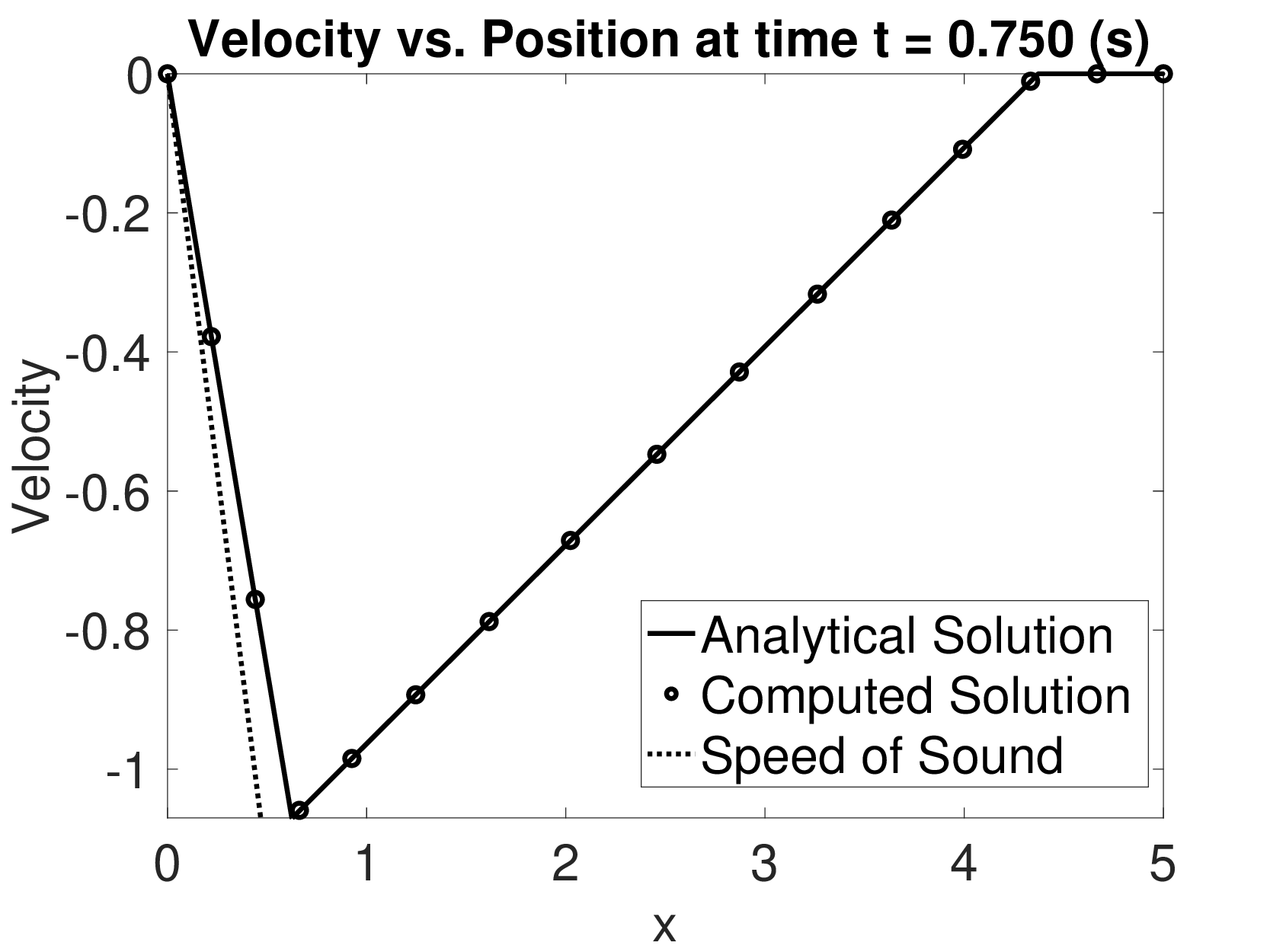}
         \caption{$t \approx 0.75$}
         \label{subfig:velocity_0.75}
     \end{subfigure}
     \hfill
     \begin{subfigure}[h]{0.49\textwidth}
         \centering
         \includegraphics[width=\textwidth]{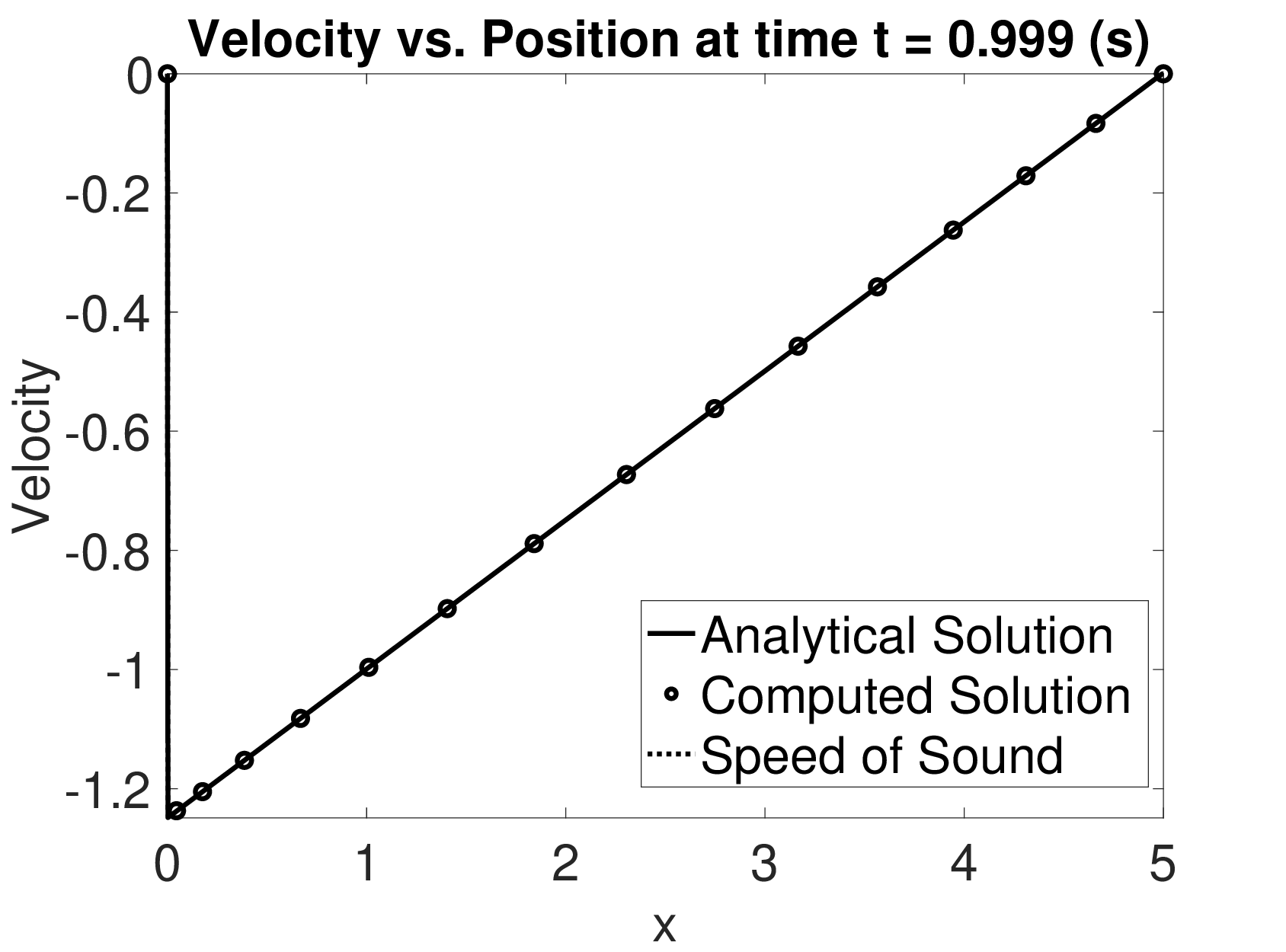}
         \caption{$t \approx 1$}
         \label{subfig:velocity_1}
     \end{subfigure}
        \caption{Analytical solution and selected representatives of the computed solution for velocity as well as the local speed of sound at approximate times $t = 0.25,0.5,0.75$ and $t = 1$.}
        \label{fig:velocity}
\end{figure}

\begin{figure}[h]
     \centering
     \begin{subfigure}[h]{0.49\textwidth}
         \centering
         \includegraphics[width=\textwidth]{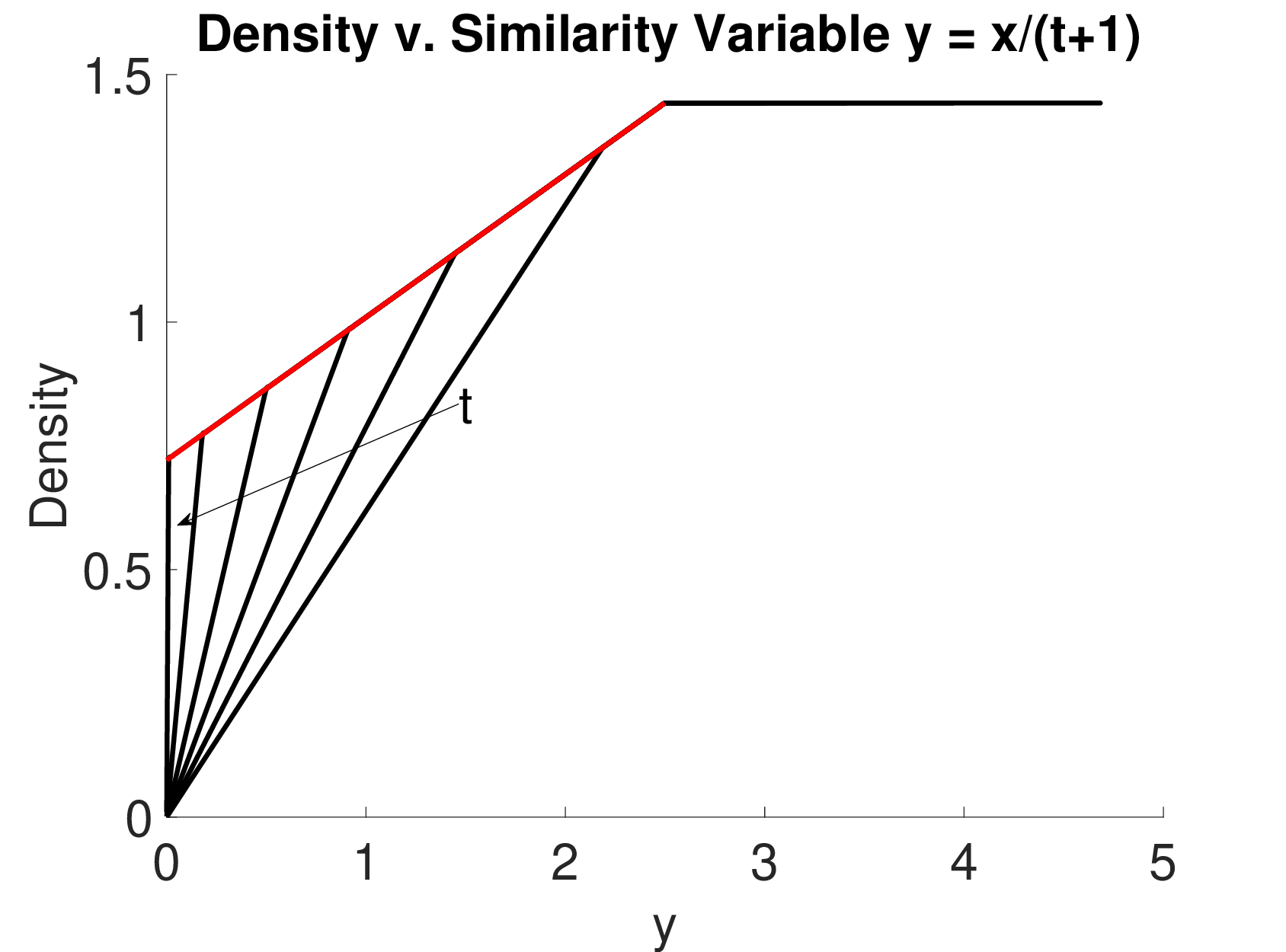}
         \caption{Density}
         \label{subfig:density_rare}
     \end{subfigure}
     \hfill
     \begin{subfigure}[h]{0.49\textwidth}
         \centering
         \includegraphics[width=\textwidth]{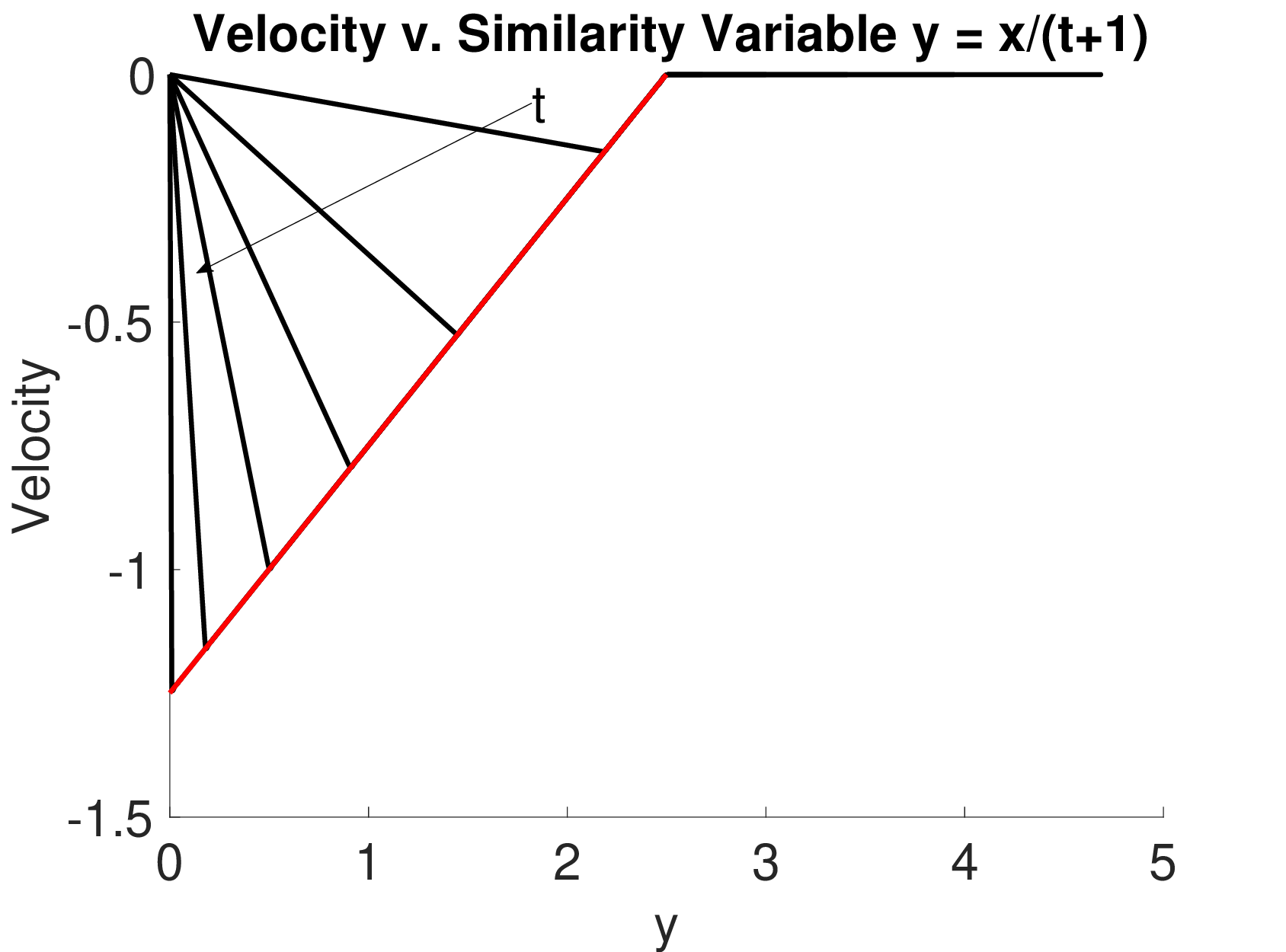}
         \caption{Velocity}
         \label{subfig:velocity_rare}
     \end{subfigure}
        \caption{Computed density and velocity plotted as a function of the similarity variable $y = \frac{x}{t+1}$ at increasing times in the direction of the arrow, as predicted in the analytical expressions for the density and velocity in equations \eqref{eq:density_cutoff_implosion} and \eqref{eq:velocity_cutoff_implosion}. The approximate times are $t = 0, 0.20, 0.33, 0.46, 0.60, 0.73, 0.87$ and $0.97$. At each time snapshot, the solution collapses onto the same middle region as highlighted in red.}
        \label{fig:rarefaction}
\end{figure}

\begin{figure}[h]
     \centering
     \begin{subfigure}[h]{0.49\textwidth}
         \centering
         \includegraphics[width=\textwidth]{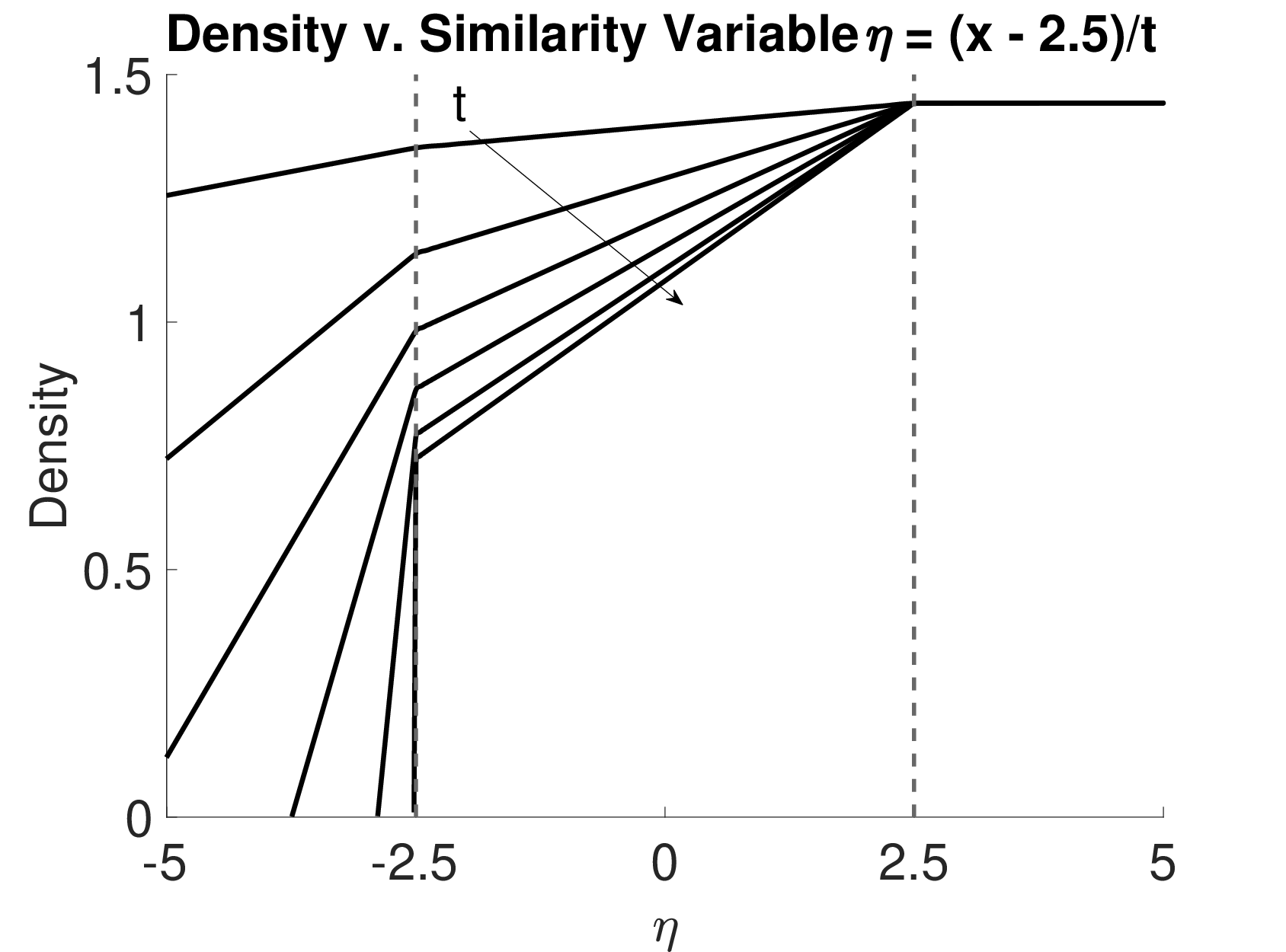}
         \caption{Density}
         \label{subfig:density_wrong_rare}
     \end{subfigure}
     \hfill
     \begin{subfigure}[h]{0.49\textwidth}
         \centering
         \includegraphics[width=\textwidth]{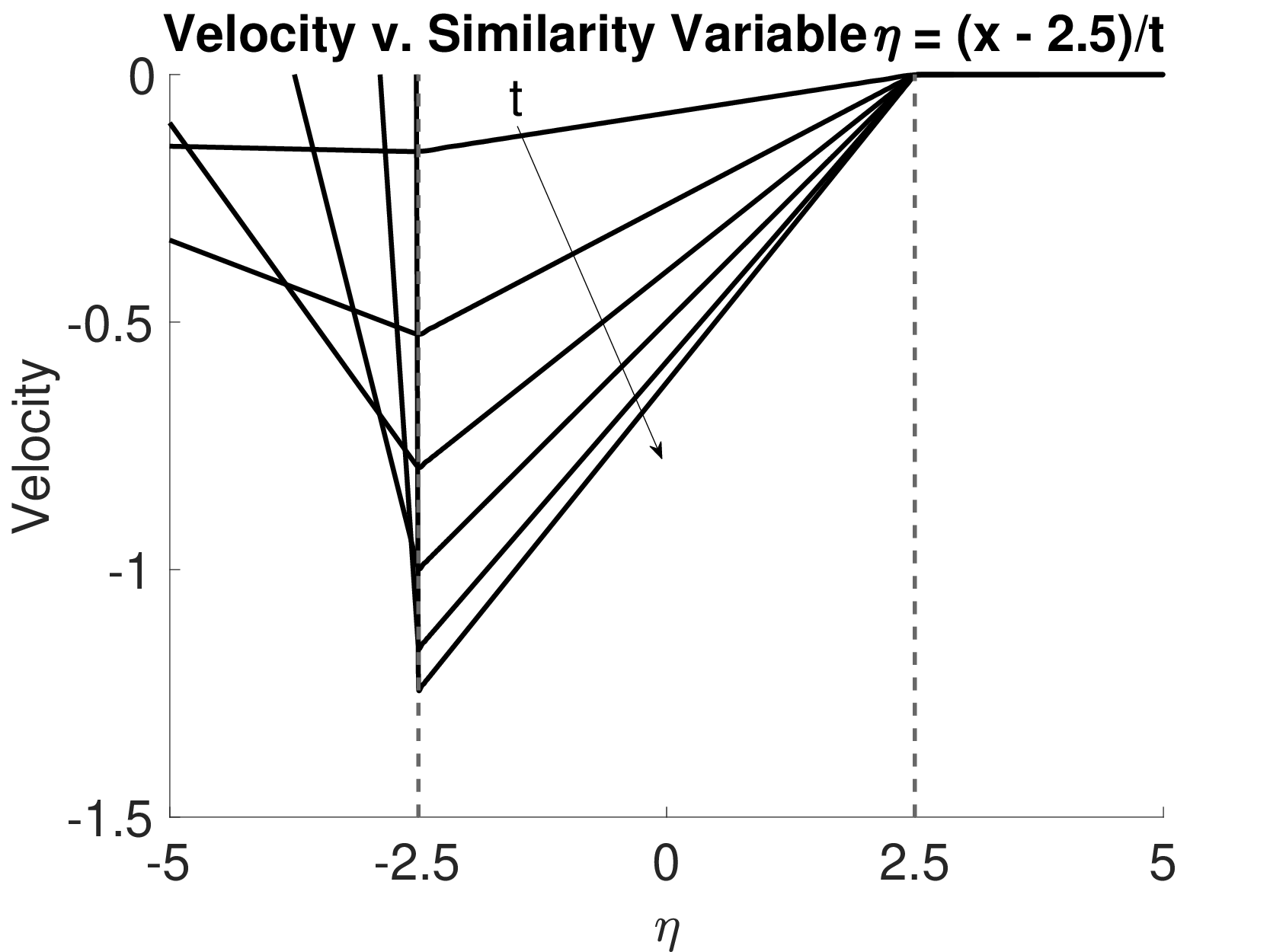}
         \caption{Velocity}
         \label{subfig:velocity_wrong_rare}
     \end{subfigure}
        \caption{Density and Velocity plotted as a function of the similarity variable $\eta = \frac{x - 2.5}{t}$ at increasing times in the direction of the arrow. The approximate times are $t = 0, 0.20, 0.33, 0.46, 0.60, 0.73, 0.87$ and $0.97$. Note the solutions at various time snapshots does not collapse onto the same curve, indicating the solution is not self-similar in $\eta$. However, the left and right boundaries of the middle region are invariant with respect to $\eta$.}
        \label{fig:rarefaction_wrong}
\end{figure}

 The initial data we use for our numerical simulations as presented in equation \eqref{eq:implosion_cutoff_initial} is shown in Figure \ref{fig:initial}. In our numerical simulations, we overlay the analytical solution for density and velocity found in equations \eqref{eq:density_cutoff_implosion} and \eqref{eq:velocity_cutoff_implosion} for comparison. We also overlay the local sound speed in the velocity plots. Figure \ref{fig:density} shows the evolution of the density at time snapshots of approximately $t = 0.25, t = 0.5, t = 0.75$, and $t = 1$. Figure \ref{fig:velocity} does the same for the velocity. In order to verify the proposed similarity variable of $y = \frac{x}{t+1}$, we plot the computed density and velocity against $y$ for various times in Figure \ref{fig:rarefaction}. Similarly, we plot the computed density and velocity against a different similarity variable $\eta = \frac{x - 2.5}{t}$ to determine the location of the rarefaction region in Figure \ref{fig:rarefaction_wrong}.

\par To quantitatively compare the numerical results to the analytical solution, we also compute the $L^1$ error norm between the two solutions. The $L^1$ norm is computed as
\begin{equation}\label{eq:L1_norm}
    \left\| y^A - y^C \right\|_1 = \frac{1}{V} \sum_{i = 1}^{N} \left| y^A(x_i) - y^C_i \right| V_i
\end{equation}
where the $A$ and $C$ superscripts represent the analytical and computed solutions respectively, $N$ is the number of nodes or volumes, $x_i$ is the location of the $i$th node or volume, $V_i$ is the zone volume, and $V$ is the sum of the zone volumes. 

%Recall velocity is defined on the nodes while density is defined on the zones, and so formula \eqref{eq:L1_norm} changes accordingly. Since the density is defined on the volumes, the position of the volumes $x_i$ are taken to be the midpoint position between nodes while the zones volumes are 

\begin{table}[h]
    \begin{tabular}{|c|c|c|}
    \hline
    $t$    & $\rho$               & $v$                  \\ \hline
    $0.25$ & $2.28 \cdot 10^{-5}$ & $2.36 \cdot 10^{-5}$ \\ \hline
    $0.50$ & $3.37 \cdot 10^{-5}$ & $3.81 \cdot 10^{-5}$ \\ \hline
    $0.75$ & $4.01 \cdot 10^{-5}$ & $5.03 \cdot 10^{-5}$ \\ \hline
    $0.99$ & $4.08 \cdot 10^{-5}$ & $6.00 \cdot 10^{-5}$ \\ \hline
    \end{tabular}
    \caption{Calculated $L^1$ error between the analytical and computed solutions for density and velocity at various times.}
    \label{tab:error}
\end{table}

\section{Discussion}\label{sec:discussion}
In this section, we will discuss the results of the numerical simulations shown in Figures \ref{fig:density}, \ref{fig:velocity}, and \ref{fig:rarefaction} and Table \ref{tab:error}.

\par First, comparing the analytical solution in equations \eqref{eq:density_cutoff_implosion} and \eqref{eq:velocity_cutoff_implosion} to the computed density in Figure \ref{fig:density} and Figure \ref{fig:velocity}, we find excellent visual agreement between the analytical and computed solutions. In particular, the computed solution captures the sharp corners corresponding to discontinuous gradient well. Diffusive numerical schemes such as Lax-Friedrichs will smoothen out such sharp edges. Sampled data points from the computed solution seem to align exactly with the analytical solution. This visual agreement is corroborated with the small computed $L^1$ error displayed in Table \ref{tab:error}. Both the qualitative and quantitative agreement instills confidence that our proposed analytical solution is correct, coupled with the discussion about unique solutions in the previous section. We next determine various features we expect to see in the numerical simulation.   

\par The computed density and velocity are split into three distinct regions. In the first region, the computed solution follows the Kidder solution. In the third region, the computed solution remain constant. In between, the computed solution connects the Kidder solution on the left to the constant solution on the right. This follows the predictions made by our derived analytical solution.

\par In our analytical solution, we posited the middle region is a rarefaction with similarity variable $y = \frac{x}{t+1}$. To verify this, we plot the computed solution at various times against this proposed similarity variable in Figure \ref{fig:rarefaction}. We see the solutions all follow a similar structure. They originate at the origin and grow linearly until collapsing onto the same curve. This continues until a later value of $y$ in which the solution becomes constant. 

\par We can explain this behavior in relation to the original structure of the computed solution. The first region represents the solutions obeying the Kidder solution. Plotted against the similarity variable, the solution follows different curves for different times. As the solution enters the rarefaction region, the solution collapses onto the same curve for every time instance. This verifies the solution is self-similar with our proposed similarity variable $y = \frac{x}{t+1}$. The third region simply indicates the solution becomes constant past a certain value.

\par Further verifying the rarefaction, we can check the endpoints of the rarefaction region by plotting the computed solution against a different similarity variable, $\eta = \frac{x - x^*}{t}$, where in our case $x^* = 2.5$. From observing the initial data, we may hypothesize the rarefaction variable should be $\eta$ as we expect a rarefaction to emerge from spatial point $x^*$ at time $t = 0$. As a stark difference from plotting the computed solutions against $y = \frac{x}{t+1}$, the solutions do not collapse onto the same curve in the rarefaction region, which supports the fact that the rarefaction similarity variable is $y$ and not $\eta$. In contrast, we do see the domain of the rarefaction is constant with respect to this variable. In particular, we see in Figure \ref{fig:rarefaction_wrong} the left endpoint of the rarefaction region corresponds to $\eta = -2.5 = -x^*$ and the right endpoint corresponds to $\eta = 2.5 = x^*$. This implies the left endpoint and right endpoint of the rarefaction region are located at $x = x^*(1-t)$ and $x = x^*(1+t)$ as we predicted.

\par The last feature we discuss is the comparison between the computed velocity and the computed sound speed. Indeed, by substituting our computed density into the formula for the local sound speed in equation \eqref{eq:sos_potential}, we observe throughout the evolution of the velocity the flow remains subsonic for all of the computation time as seen in Figure \ref{fig:velocity}. 

%Superflous paragraph
%\par By comparing our analytical solution to numerical simulations of the potential flow equations with cutoff implosion initial data, we verify that our proposed analytical solution is correct and important structures of our analytical solution - such as the proposed similarity variable - are also observed in the numerically computed solutions as well.  

\section{Conclusion}\label{sec:conclusion}
\par In this work, we construct an explicit solution to the cutoff implosion problem. We start by showing how to derive the generalization of the homogeneous implosion solutions of \cite{Kidder1974}. By specifically considering the one dimensional case, we define the cutoff implosion initial data and show that an exact solution can be constructed with a rarefaction expanding outward from the initial cutoff point in space. We observe distinguished features of the solution, such as its nonstandard rarefaction variable and the suppression of the imploding solution. To compare with our analytic solution, we numerically simulate the cutoff implosion problem using a Lagrangian leapfrog scheme, showing excellent agreement.

\par The strategy of cutting off the initial data of the Kidder problem with a constant is physically motivated. The Kidder solution is unbounded at infinity which is unrealistic. Moreover the original solution blows up at all points in space as time approaches the implosion time.  Instead, we cut off the pressure and density in the far field with a constant.
This corresponds physically to instead submerging the gas in an ambient environment. 

\par However we discover that the effect of the cutoff leads to a growing rarefaction wave that emerges from the cutoff position and eventually dominates
the imploding part of the solution. This is an interesting phenomenon on its own. The rarefaction similarity variable arising in the cutoff implosion solution is nonstandard. Our rarefaction variable also encodes the original implosion time within itself. In the derivation of the rarefaction, this arises as a byproduct of taking the left state to be the time varying Kidder solution. The emergence of the unusual rarefaction variable in a fairly simple problem is interesting and warrants further exploration. We began the investigation by hypothesizing whether the finite speed of propagation of the cutoff will reach the origin and disturb the Kidder solution. However, regardless of where the cutoff is positioned, the Kidder solution vanishes.

\par The theoretical work done in this paper is supported with numerical experiments. Such numerical simulations can be also performed for the imploding solutions of \cite{Kidder1974}, suggesting the imploding solutions of \cite{Kidder1974} are numerically stable - unlike the imploding solutions investigated in \cite{biasi_compute}. The numerical stability of the cutoff implosion solution may imply its nonlinear stability as well, which in turn may give more insight into the physical existence of the solution. Numerical investigation into imploding solutions of this type could be an avenue of further discovery into implosions, similar to the work done in this paper. One such exploratory avenue is the cutoff implosion problem for higher dimensions.

\par The 1D cutoff implosion problem may be a good numerical test problem. Some
features making the cutoff implosion problem a desirable numerical test problem are the similarity solutions present in the solution before $t = 1$. A candidate numerical code should be able to identify the similarity variable $y$ within. Additionally, a candidate numerical code can be evaluated on how it resolves the solution past the
shock that forms at $t = 1$, of which we know the theoretical intermediate state and shock speed. The multifaceted dynamics of the cutoff implosion problem make it a rich test for potential numerical codes. We also note the proposed post shock solution is quite similar to the Noh problem \cite{NohArtificialViscosity}, a common numerical test problem. In particular, the zero velocity and higher density in the shocked region is seen in both the cutoff implosion problem and the Noh problem.

%\par In higher dimensions, it is currently unclear whether the constant cutoff also suppresses the implosion. The argument showing how the rarefaction dominates the implosion in this work does not extend neatly to higher dimensions, as the presence of the source term in the radially symmetric, higher dimension Euler equations prevents the typical rarefaction analysis for systems of hyperbolic conservation laws used in one dimension. Such rarefaction analysis gives the closed form solution for the rarefaction in this problem. The closed form is integral in showing the suppression of the implosion. We foresee any investigation into the higher dimension problem to require a different approach than the rarefaction analysis done in this work. Studying the constant cutoff in higher dimensions is important though as it may give richer insight to the physical understanding of the Kidder solution.

\section*{Acknowledgments}
The authors would like to thank Alexis Vasseur for his insights in the application of weak-strong uniqueness results. The authors would like to thank the reviewer's suggestion to use the additional entropy pair available for $\gamma = 3$.
%\section{Appendix}
%\appendix
%\section{Linear Velocity Derivation}
%\section{Hyperbolic Conservation Law Theory}
\bibliographystyle{plain}
\bibliography{sources}
\end{document}